\documentclass[11pt]{article}
\usepackage{hyperref, amsmath, amssymb, amsthm, fullpage, enumitem, mathtools, microtype, tikz, bm, subcaption, alphalph}
\usepackage[noabbrev, capitalize]{cleveref}

\newtheorem{theorem}{Theorem}[section]
\newtheorem{corollary}[theorem]{Corollary}
\newtheorem{lemma}[theorem]{Lemma}
\newtheorem{proposition}[theorem]{Proposition}
\newtheorem{problem}[theorem]{Problem}

\theoremstyle{definition}
\newtheorem{definition}[theorem]{Definition}

\AddToHook{env/corollary/begin}{\crefalias{theorem}{corollary}}
\AddToHook{env/lemma/begin}{\crefalias{theorem}{lemma}}
\AddToHook{env/proposition/begin}{\crefalias{theorem}{proposition}}
\AddToHook{env/problem/begin}{\crefalias{theorem}{problem}}
\AddToHook{env/definition/begin}{\crefalias{theorem}{definition}}

\theoremstyle{remark}

\newtheorem*{example}{Example}
\newtheorem*{claim*}{Claim}

\crefname{claim}{Claim}{Claims}

\newlist{kase}{enumerate}{4} 
\setlistdepth{4}
\setlist[kase,1]{label=(\arabic*)}
\setlist[kase,2]{label=(\arabic{kasei}.\arabic*)}
\setlist[kase,3]{label=(\arabic{kasei}.\arabic{kaseii}.\arabic*)}
\setlist[kase,4]{label=(\arabic{kasei}.\arabic{kaseii}.\arabic{kaseiii}.\arabic*)}

\newcommand{\figurehc}[2][1]{
\begin{center}
    \centering
    \begin{tikzpicture}[scale=#1]
        #2
    \end{tikzpicture}
\end{center}}

\newcommand{\gfour}[2]{%
  \begin{subfigure}{2cm}
    \centering
    \begin{tikzpicture}
      #2
    \end{tikzpicture}
    \caption{}\label{fig:g4-#1}
  \end{subfigure}%
}

\newcommand{\beineke}[2]{%
  \begin{subfigure}{3cm}
    \centering
    \begin{tikzpicture}[baseline=(0.base)]
      #2
    \end{tikzpicture}
    \caption{}\label{fig:beineke-#1}
  \end{subfigure}%
}

\newcommand{\dset}[2]{\left\{{#1}\colon{#2}\right\}}
\newcommand{\sset}[1]{\left\{{#1}\right\}}
\usepackage{graphicx}
\newcommand{\upbowtie}{\mathrel{\rotatebox{90}{\scalebox{0.75}{$\bm{\bowtie}$}}}}
\newcommand{\ol}[1]{\overline{\vphantom{b}#1}}

\newcommand{\N}{\mathbb{N}}
\newcommand{\D}{\mathcal{D}}

\newcommand{\gmn}{\mathsf{GM}_n}
\newcommand{\hjn}{\mathsf{HJ}_n}
\newcommand{\hjone}{\mathsf{HJ}_1^+}
\newcommand{\iddtwo}{I(\D)_{\D_{\sharp 2}}}

\usetikzlibrary{decorations.pathreplacing, calc}
\definecolor{litegray}{RGB}{220,220,220}
\tikzstyle{vertex}=[circle, draw, fill=white, inner sep=0pt, minimum width=4pt]
\tikzstyle{root-vertex}=[circle, draw, fill=darkgray, inner sep=0pt, minimum width=4pt]
\tikzset{
  split-vertex/.style={
    circle, draw,
    inner sep=0pt,
    minimum size=4pt,
    fill=darkgray,
    path picture={
      \begin{scope}
        \clip (path picture bounding box.center) circle[radius=2pt];
        \fill[white]
          (path picture bounding box.south west) --
          (path picture bounding box.north west) --
          (path picture bounding box.north east) -- cycle;
      \end{scope}
    }
  }
}

\tikzstyle{braket}=[decorate,decoration={brace,amplitude=10pt},xshift=0pt,yshift=-10pt, black]
\tikzset{every picture/.append style={scale=.8, thick}}
\tikzstyle{f-vertex}=[circle, draw=black, fill=white, text=black, font=\footnotesize, inner sep=0pt, minimum width=12pt]
\newcommand{\drawlabeledges}[1]{\foreach \u/\v/\l in {#1} {\draw (\u) -- (\v) node[midway, fill=white, inner sep=2pt, font=\footnotesize] {$\l$};}}
\newcommand{\drawedges}[1]{\foreach \u/\v in {#1} {\draw (\u) -- (\v);}}
\newcommand{\defxy}[1]{\foreach \u/\x/\y in {#1} {\coordinate (\u) at (\x,\y);}}
\newcommand{\defpolar}[1]{\foreach \u/\a/\r in {#1} {\coordinate (\u) at (\a:\r);}}
\newcommand{\nodesvertex}[1]{\foreach \u in {#1} {\node[vertex] at (\u) {};}}
\newcommand{\nodesrootvertex}[1]{\foreach \u in {#1} {\node[root-vertex] at (\u) {};}}
\newcommand{\nodessplitvertex}[1]{\foreach \u in {#1} {\node[split-vertex] at (\u) {};}}
\newcommand{\nodeslabelvertex}[1]{\foreach \u/\l in {#1} {\node[f-vertex] at (\u) {$\l$};}}
\newcommand{\nodesvertexbadge}[1]{\foreach \u/\l in {#1} {\node[rectangle, fill=litegray, rounded corners, font=\footnotesize, inner sep=3pt, minimum height=1.4em] at (\u) {$\l$};}}
\newcommand{\nodesvertexbadgehollow}[1]{\foreach \u/\l in {#1} {\node[rectangle, fill=white, rounded corners, font=\footnotesize, inner sep=3pt] at (\u) {$\l$};}}
\newcommand{\drawedgesbadge}[1]{\foreach \u/\v/\l in {#1} {\draw (\u) -- (\v) node[midway, fill=litegray, rounded corners, font=\footnotesize, inner sep=3pt, minimum height=1.4em] {$\l$};}}

\begin{document}

\title{Subcubic graphs without eigenvalues in $(-1, 1)$}
\author{
    Shenwei Huang\thanks{School of Mathematical Sciences and LPMC, Nankai University, Tianjin 300071, P.R.~China. Email: {\tt shenweihuang@nankai.edu.cn}.}
    \and Zilin Jiang\thanks{School of Mathematical and Statistical Sciences, and School of Computing and Augmented Intelligence, Arizona State University, Tempe, AZ 85281, USA. Email: {\tt zilinj@asu.edu}. Supported in part by the Simons Foundation through its Travel Support for Mathematicians program and by U.S.\ taxpayers through NSF grant 2451581.}
}
\date{}

\maketitle

\begin{abstract}
    Guo and Royle recently classified the connected cubic graphs without eigenvalues of their adjacency matrix in the open interval $(-1, 1)$, and raised the question of extending their classification to graphs of maximum degree at most $3$. Together with their cubic classification, our result fully answers this question by characterizing all connected subcubic graphs that are not cubic and have no eigenvalues in $(-1,1)$. We show that exactly two infinite families and seven sporadic examples occur, and that every sporadic graph has at most $18$ vertices.

    To obtain this complete classification, we build a bridge between spectral graph theory and structural graph theory for graphs whose adjacency matrices, after selected diagonal entries are changed to $-1$, have smallest eigenvalue at least $-2$. This generalizes the classical theorem of Cameron, Goethals, Seidel and Shult for graphs with smallest eigenvalue at least $-2$.

    As a consequence, we prove that $(-1,1)$ is a maximal spectral gap set for the class of connected subcubic graphs. Guo and Royle, answering a question of Koll\'ar and Sarnak, established this maximality for connected cubic graphs. Our result generalizes their conclusion to the subcubic setting.
\end{abstract}

\section{Introduction} \label{sec:intro}

Guo and Royle~\cite{GR26} recently classified the connected cubic graphs, that is, $3$-regular graphs, without eigenvalues in $(-1, 1)$. Here and throughout, all graph eigenvalues are those of the adjacency matrix. They found two infinite families, the Koll\'ar--Sarnak graphs and the Guo--Mohar graphs, together with $13$ sporadic examples, including the Petersen graph, the Heawood graph, the M\"obius--Kantor graph, and the Desargues graph (see \cref{fig:sporadic-cubic} for four of the sporadic graphs).

\begin{figure}[t]
    \centering
    \begin{tikzpicture}[scale=2]
        \foreach \i in {0,...,4} {
            \node[vertex] (A\i) at ({360/5 * \i+18}:1) {};
        }
        \foreach \i in {0,...,4} {
            \pgfmathtruncatemacro{\next}{Mod(\i + 1, 5)}
            \draw (A\i) -- (A\next);
        }
        \foreach \i in {0,...,4} {
            \node[vertex] (B\i) at ({360/5 * \i+18}:0.5) {};
            \draw (A\i) -- (B\i);
        }
        \foreach \i in {0,...,4} {
            \pgfmathtruncatemacro{\next}{Mod(\i + 2, 5)}
            \draw (B\i) -- (B\next);
        }
    \end{tikzpicture}\qquad
    \begin{tikzpicture}[scale=2]
        \foreach \i in {0,...,13} {
            \node[vertex] (\i) at ({360/14 * \i}:1) {};
        }
        \foreach \i in {0,...,13} {
            \pgfmathtruncatemacro{\next}{Mod(\i + 1, 14)}
            \draw (\i) -- (\next);
        }
        \foreach \i in {0,...,13} {
            \pgfmathtruncatemacro{\iseven}{Mod(\i, 2)}
            \ifnum\iseven=0
                \pgfmathtruncatemacro{\target}{Mod(\i + 5, 14)}
            \else
                \pgfmathtruncatemacro{\target}{Mod(\i - 5, 14)}
            \fi
            \draw (\i) -- (\target);
        }
    \end{tikzpicture}\qquad
    \begin{tikzpicture}[scale=2]
        \foreach \i in {0,...,7} {
            \node[vertex] (A\i) at ({360/8 * \i+22.5}:1) {};
        }
        \foreach \i in {0,...,7} {
            \pgfmathtruncatemacro{\next}{Mod(\i + 1, 8)}
            \draw (A\i) -- (A\next);
        }
        \foreach \i in {0,...,7} {
            \node[vertex] (B\i) at ({360/8 * \i+22.5}:0.5) {};
            \draw (A\i) -- (B\i);
        }
        \foreach \i in {0,...,7} {
            \pgfmathtruncatemacro{\next}{Mod(\i + 3, 8)}
            \draw (B\i) -- (B\next);
        }
    \end{tikzpicture}\qquad
    \begin{tikzpicture}[scale=2]
        \foreach \i in {0,...,9} {
            \node[vertex] (A\i) at ({360/10 * \i}:1) {};
        }
        \foreach \i in {0,...,9} {
            \pgfmathtruncatemacro{\next}{Mod(\i + 1, 10)}
            \draw (A\i) -- (A\next);
        }
        \foreach \i in {0,...,9} {
            \node[vertex] (B\i) at ({360/10 * \i}:0.5) {};
            \draw (A\i) -- (B\i);
        }
        \foreach \i in {0,...,9} {
            \pgfmathtruncatemacro{\next}{Mod(\i + 3, 10)}
            \draw (B\i) -- (B\next);
        }
    \end{tikzpicture}
    \caption{Petersen graph, Heawood graph, M\"obius--Kantor graph, and Desargues graph} \label{fig:sporadic-cubic}
\end{figure}

The two infinite families are built from the \emph{twisted ladder} $\upbowtie_n$ on $n$ rungs shown in \cref{fig:twisted-ladder}. For $n \ge 1$, the Koll\'ar--Sarnak graph $\mathsf{KS}_n$ is obtained from the twisted ladder $\upbowtie_n$ by adding edges $ab$ and $yz$, and, for $n \ge 2$, the Guo--Mohar graph $\mathsf{GM}_n$ is obtained from $\upbowtie_n$ by adding $ay$ and $bz$ instead.

\begin{figure}[b]
    \centering
    \begin{tikzpicture}
        \defxy{a1/0/0, b1/0/1, c1/1/0, d1/1/1, a2/2/0, b2/2/1, c2/3/0, d2/3/1, a3/4/0, b3/4/1, c3/6/0, d3/6/1, a4/7/0, b4/7/1, c4/8/0, d4/8/1, a5/9/0, b5/9/1, c5/10/0, d5/10/1}
        \drawedges{a1/c1, c1/b1, b1/d1, d1/a1, a2/c2, c2/b2, b2/d2, d2/a2, a4/c4, c4/b4, b4/d4, d4/a4, a5/c5, c5/b5, b5/d5, d5/a5, c1/a2, d1/b2, c2/a3, d2/b3, c3/a4, d3/b4, c4/a5, d4/b5}
        \draw[dashed] (a3)--(c3);
        \draw[dashed] (b3)--(d3);
        \node[vertex, label=left:$b$] at (a1) {};
        \node[vertex, label=left:$a$] at (b1) {};
        \nodesvertex{c1,d1,a2,b2,c2,d2,a3,b3,c3,d3,a4,b4,c4,d4,a5,b5}
        \node[vertex, label=right:$z$] at (c5) {};
        \node[vertex, label=right:$y$] at (d5) {};
        \draw [braket] (0,1.5) -- (10,1.5) node [black,midway,yshift=17pt] {\footnotesize $2n$};
    \end{tikzpicture}
    \caption{Twisted ladder $\upbowtie_n$ on $n$ rungs, which has $4n$ vertices.} \label{fig:twisted-ladder}
\end{figure}

Guo and Royle asked for a classification of connected subcubic graphs, that is, connected graphs of maximum degree at most $3$, without eigenvalues in $(-1,1)$. Their cubic classification relies essentially on regularity, and they observed that their techniques do not directly extend to the subcubic setting.

Our main result answers their question in the complementary setting of connected subcubic graphs that are not cubic. The passage from cubic to subcubic graphs is not merely a matter of treating leftover cases: the loss of regularity forces one to retain extra structural information about the vertices of degree $2$. While \cite{GR26} uses computer-assisted methods, our proof is purely theoretical.

\begin{theorem} \label{thm:main}
    The family of connected subcubic graphs that are not cubic and have no eigenvalues in $(-1,1)$ consists of the bipartite graphs $\hjn$ for $n \ge 1$ in \cref{fig:bipartite-girth-4}, the non-bipartite graphs $\hjn'$ for $n \ge 1$ in \cref{fig:non-bipartite}, and seven sporadic examples --- the single-edge graph $K_2$, the triangle graph $K_3$, the graph in \cref{fig:bipartite-girth-4-sporadic}, and the four graphs in \cref{fig:bipartite-girth-6}.
\end{theorem}

This classification is technically different from the cubic case. In the cubic setting, regularity makes the spectral condition uniform across all vertices. In the subcubic setting, vertices of degree $2$ must be recorded separately, and this additional information is used throughout the proof.

\begin{figure}
    \centering
    \begin{tikzpicture}
        \defxy{a1/0/0,b1/0/1,c1/1/0,d1/1/1,a2/2/0,b2/2/1,c2/3/0,d2/3/1,a3/4/0,b3/4/1,c3/6/0,d3/6/1,a4/7/0,b4/7/1,c4/8/0,d4/8/1,a5/9/0,b5/9/1,c5/10/0,d5/10/1,b0/-2/0,c0/-1/0,d0/-1/1,a6/11/0,b6/11/1,c6/12/0}
        \drawedges{b1/d0,d0/b0,b0/c0,c0/a1,a1/c1,c1/b1,b1/d1,d1/a1,a2/c2,c2/b2,b2/d2,d2/a2,a4/c4,c4/b4,b4/d4,d4/a4,c5/b5,b5/d5,d5/a5,a5/c5,c5/a6,a6/c6,c6/b6,b6/d5,c1/a2,d1/b2,c2/a3,d2/b3,c3/a4,d3/b4,c4/a5,d4/b5}
        \draw[dashed] (a3)--(c3);
        \draw[dashed] (b3)--(d3);
        \nodesvertex{a1,b1,c1,d1,a2,b2,c2,d2,a3,b3,c3,d3,a4,b4,c4,d4,a5,b5,c5,d5,c0,d0,b0,a6,b6,c6}
        \draw [braket] (-1,1.5) -- (11,1.5) node [black,midway,yshift=17pt] {\footnotesize $2n+2$};
    \end{tikzpicture}
    \caption{Bipartite graphs $\hjn$ without eigenvalues in $(-1,1)$.}
    \label{fig:bipartite-girth-4}
\end{figure}

\begin{figure}[b]
    \centering
    \begin{tikzpicture}
        \defxy{b0/-2/0,c0/-1/0,d0/-1/1,a1/0/0,b1/0/1,c1/1/0,d1/1/1,a2/2/0,b2/2/1,c2/3/0,d2/3/1,a3/4/0,b3/4/1,c3/6/0,d3/6/1,a4/7/0,b4/7/1,c4/8/0,d4/8/1}
        \drawedges{b1/d0,d0/b0,b0/c0,c0/a1,a1/c1,c1/b1,b1/d1,d1/a1,a2/c2,c2/b2,b2/d2,d2/a2,a4/c4,c4/b4,b4/d4,d4/a4,c1/a2,d1/b2,c2/a3,d2/b3,c3/a4,d3/b4,c4/d4}
        \draw[dashed] (a3)--(c3);
        \draw[dashed] (b3)--(d3);
        \nodesvertex{a1,b1,c1,d1,a2,b2,c2,d2,a3,b3,c3,d3,a4,b4,c4,d4,b0,c0,d0}
        \draw [braket] (-1,1.5) -- (8,1.5) node [black,midway,yshift=17pt] {\footnotesize $2n+1$};
    \end{tikzpicture}
    \caption{Non-bipartite graphs $\hjn'$ without eigenvalues in $(-1,1)$.}
    \label{fig:non-bipartite}
\end{figure}

We now explain the motivation for this classification in terms of \emph{spectral gap sets}. Let $\mathcal{C}$ be a class of graphs, and let $a$ and $b$ denote the infimum and supremum, respectively, of all eigenvalues of graphs in $\mathcal{C}$. An open subset $I$ of $[a,b]$ is called a \emph{spectral gap set} for $\mathcal{C}$ if there exist infinitely many graphs in $\mathcal{C}$ whose eigenvalues completely avoid $I$. A spectral gap set $I$ for $\mathcal{C}$ is said to be \emph{maximal} if it is not properly contained in any other spectral gap set for $\mathcal{C}$.

The notions of spectral gap sets and maximality were introduced by Koll\'ar and Sarnak in \cite{KS21}. In the case where $\mathcal{C}$ is the class of connected cubic graphs, they pointed out that several classical spectral graph-theoretic phenomena can be reformulated as statements about spectral gap sets. For instance, the classical fact that line graphs have smallest eigenvalue at least $-2$, together with the observation that the line graph of a subdivision of a cubic graph is again cubic, implies that $[-3,-2)$ is a spectral gap set for $\mathcal{C}$. Likewise, Chiu's construction of cubic Ramanujan graphs \cite{C92} implies that $[-3,-2\sqrt2) \cup (2\sqrt2, 3)$ is a spectral gap set.

Koll\'ar and Sarnak \cite{KS21}, and earlier Guo and Mohar \cite{GM14}, showed that the interior interval $(-1,1)$ is also a spectral gap set for the class of connected cubic graphs. Subsequently, Guo and Royle \cite{GR26} proved that this gap is maximal by classifying all connected cubic graphs with no eigenvalues in $(-1,1)$. Our main theorem extends this picture beyond the cubic setting: combining \cref{thm:main} with the classification of Guo and Royle yields a complete description of all connected subcubic graphs with no eigenvalues in $(-1,1)$. As a consequence, we obtain the following strengthening of their result on spectral gap sets.

\begin{corollary} \label{cor:gap}
    The interval $(-1,1)$ is a maximal spectral gap set for the class of connected subcubic graphs.
\end{corollary}

\section{Proof ideas} \label{sec:idea}

We begin by disposing of the trivial low-degree cases.

\begin{proposition} \label{thm:min-deg-1}
    For every connected graph without eigenvalues in $(-1,1)$, if its minimum degree is at most $1$, then it is isomorphic to the single edge graph $K_2$.
\end{proposition}

\begin{proof}
    Clearly, the trivial graph $K_1$ has an eigenvalue in $(-1,1)$. Suppose that $H$ is a connected graph without eigenvalues in $(-1,1)$, and suppose that $H$ has a pendant edge $uv$ such that $\deg_H(u) = 1$. Assume for the sake of contradiction that $H$ is not isomorphic to $K_2$, that is, $\deg_H(v) \ge 2$. Let $w$ be another neighbor of $v$. Since $H$ has no eigenvalues in $(-1,1)$, the matrix $A_H^2 - I$ is positive semidefinite. In particular, its principal submatrix
    \[
        \begin{pmatrix}
            0 & 1 \\
            1 & \deg_H(w)-1
        \end{pmatrix}
    \]
    induced by $\sset{u, w}$ is positive semidefinite. However, its determinant is $-1$.
\end{proof}

We next introduce the language used to encode walks of length $2$ while retaining the vertices of degree $2$.

\begin{definition}[Distance-two multigraph]
    The \emph{distance-two multigraph} of a graph $H$ is the multigraph with vertex set $V(H)$ in which the multiplicity of the edge between two vertices equals the number of distinct paths of length $2$ between them in $H$.
\end{definition}

Since we must now allow vertices of degree $2$ in $H$, it becomes necessary to keep track of such vertices separately, as they affect the diagonal entries of $A_H^2$. This leads us to work with rooted multigraphs.

\begin{definition}[Rooted multigraph]
    A \emph{rooted multigraph} $G_R$ is a multigraph $G$ equipped with a distinguished subset $R \subseteq V(G)$ of vertices, called the \emph{roots} (depicted by solid circles).
\end{definition}

\begin{definition}[Rooted distance-two induced subgraph]
    Given a graph $H$ with minimum degree at least $2$ and a vertex subset $U \subseteq V(H)$, the \emph{rooted distance-two subgraph} $G_R$ of $H$ induced by $U$ is a subgraph $G$ of the distance-two multigraph of $H$ induced by $U$, where the root set $R$ consists of those vertices in $U$ that have degree $2$ in $H$. A \emph{rooted distance-two component} of $H$ is a rooted distance-two subgraph of $H$ induced by a connected component of the distance-two multigraph of $H$.
\end{definition}

We now express $A_H^2 - I \succeq 0$ in terms of these rooted distance-two induced subgraphs. Recall that the adjacency matrix of a multigraph records edge multiplicities.

\begin{definition}[Modified adjacency matrix]
    Given a rooted multigraph $G_R$, the \emph{modified adjacency matrix} of $G_R$ is the adjacency matrix of $G$, except that the diagonal entries corresponding to the vertices in $R$ are decreased from $0$ to $-1$.
\end{definition}

\begin{lemma} \label{lem:psd}
    For every subcubic graph $H$ with minimum degree $2$, the graph $H$ has no eigenvalues in $(-1,1)$ if and only if, for every rooted distance-two induced subgraph $G_R$ of $H$, the modified adjacency matrix of $G_R$ has smallest eigenvalue at least $-2$.
\end{lemma}

\begin{proof}
    Clearly, $H$ has no eigenvalues in $(-1,1)$ if and only if $A_H^2 - I \succeq 0$ if and only if every principal submatrix of $A_H^2 - I$ is positive semidefinite. Finally, notice that, for every vertex subset $U$ of $H$, the principal submatrix of $A_H^2 - 3I$ induced by $U$ is equal to the modified adjacency matrix of the rooted distance-two subgraph of $H$ induced by $U$.
\end{proof}

The preceding lemma is the point at which the connection with the work of Guo and Royle becomes precise. In the cubic case, the condition $A_H^2 \succeq I$ implies that $A_H^2-3I$ is the adjacency matrix of a multigraph whose smallest eigenvalue is at least $-2$, allowing one to invoke the classification theorem of Cameron, Goethals, Seidel and Shult \cite{CGSS76}.

For subcubic graphs, this direct translation is no longer available: vertices of degree $2$ contribute diagonal entries equal to $-1$ in the modified adjacency matrices above. Thus one needs structural information about graphs whose modified adjacency matrices have smallest eigenvalue at least $-2$. This provides a new bridge between spectral and structural graph theory, generalizing the classical Cameron--Goethals--Seidel--Shult bridge for ordinary adjacency matrices with no modified diagonal entries.

We shall use \cref{lem:psd} and the following computation to forbid certain rooted distance-two induced subgraphs.

\begin{proposition} \label{lem:forb-bipartite}
    None of the rooted multigraphs in \cref{fig:forb-girth-4} has a modified adjacency matrix whose smallest eigenvalue is at least $-2$.
\end{proposition}

\begin{figure}
    \centering
    \gfour{a}{
        \defxy{a/0/0,b/0/1}
        \draw (a) to[bend right=30] (b) to[bend right=30] (a);
        \nodessplitvertex{b}
        \nodesrootvertex{a}
    }
    \gfour{b}{
        \defxy{a/0/0,b/0/1}
        \draw (a) -- (b);
        \draw (a) to[bend right=45] (b) to[bend right=45] (a);
        \nodesvertex{a,b}
    }
    \gfour{c}{
        \defxy{a/0/0,b/0/1,c/1/1}
        \draw (a) to[bend right=30] (b) to[bend right=30] (a);
        \draw (b) -- (c);
        \nodessplitvertex{a,c}
        \nodesvertex{b}
    }
    \gfour{d}{
        \defxy{a/0/0,b/0/1,c/1/1}
        \draw (a) -- (c);
        \draw (a) to[bend right=30] (b) to[bend right=30] (a);
        \draw (b) to[bend right=30] (c) to[bend right=30] (b);
        \nodessplitvertex{a,b,c}
    }
    \gfour{e}{
        \defxy{a/0/0,b/0/1,c/1/1,d/1/0}
        \draw (a) to (b) to (c) to (d);
        \nodesvertex{c}
        \nodesrootvertex{b}
        \nodessplitvertex{a,d}
    }
    \gfour{f}{
        \defxy{a/0/0,b/0/1,c/1/1,d/1/0}
        \draw (a) to (b) to (c) to (d) to (a);
        \nodesvertex{c,d}
        \nodesrootvertex{b}
        \nodessplitvertex{a}
    }
    \gfour{g}{
        \defxy{a/0/1,b/0/0,c/1/1}
        \draw (b) to (a) to (c);
        \nodesrootvertex{a,b,c}
    }
    \caption{Rooted multigraphs whose modified adjacency matrices have smallest eigenvalue less than $-2$. Solid circles indicate roots; half-filled circles may or may not be roots; hollow circles are non-roots.} \label{fig:forb-girth-4}
\end{figure}

\begin{proof}
    When two rooted multigraphs $G_R$ and $G_S$ only differ in their root sets $R$ and $S$, if $R \subseteq S$, then the modified adjacency matrices $A$ and $B$ of $G_R$ and $G_S$ satisfy $A \succeq B$. We may assume that the half-filled circles in \cref{fig:forb-girth-4} are non-roots. One can then compute the determinants of the modified adjacency matrices of the rooted graphs plus $2I$:
    \begin{gather*}
        \det\begin{pmatrix}
            2 & 2 \\
            2 & 1
        \end{pmatrix} = -2,\quad
        \det\begin{pmatrix}
            2 & 3 \\
            3 & 2
        \end{pmatrix} = -5,\quad
        \det\begin{pmatrix}
            2 & 2 & 1 \\
            2 & 2 & 0 \\
            1 & 0 & 2
        \end{pmatrix} = -2,\quad
        \det\begin{pmatrix}
            2 & 2 & 2 \\
            2 & 2 & 1 \\
            2 & 1 & 2
        \end{pmatrix} = -2,\\
        \det\begin{pmatrix}
            1 & 1 & 1 & 0 \\
            1 & 2 & 0 & 0 \\
            1 & 0 & 2 & 1 \\
            0 & 0 & 1 & 2
        \end{pmatrix} = -1,\quad
        \det\begin{pmatrix}
            1 & 1 & 1 & 0 \\
            1 & 2 & 0 & 1 \\
            1 & 0 & 2 & 1 \\
            0 & 1 & 1 & 2
        \end{pmatrix} = -4,\quad
        \det\begin{pmatrix}
            1 & 1 & 0 \\
            1 & 1 & 1 \\
            0 & 1 & 1
        \end{pmatrix} = -1,
    \end{gather*}
    all of which are negative.
\end{proof}

The preceding spectral preparation allows us to focus on the following purely structural classification of bipartite graphs of girth $4$.

\begin{theorem} \label{thm:bipartite-girth-4}
    For every connected subcubic bipartite graph $H$ of minimum degree $2$ and girth $4$, if none of the rooted multigraphs in \crefrange{fig:g4-a}{fig:g4-f} is a rooted distance-two induced subgraph of $H$, then $H$ is isomorphic to $\mathsf{HJ}_n$ for some $n \ge 1$, or $\hjone$ defined in \cref{fig:bipartite-girth-4-sporadic}.
\end{theorem}

For bipartite graphs $H$ of girth at least $6$, the distance-two multigraph is in fact a simple graph. In the presence of vertices of degree $2$, although \cref{lem:psd} implies that each distance-two component $G_R$ of $H$ has smallest eigenvalue at least $-2$, this implication does not exploit the fact that some diagonal entries of the modified adjacency matrix of $G_R$ can be equal to $-1$. The following lemma provides the key bridge that allows us to apply the classical classification theorem of Cameron, Goethals, Seidel, and Shult~\cite{CGSS76} while retaining the information carried by the $-1$ diagonal entries.

\begin{figure}
    \centering
    \begin{tikzpicture}
        \defxy{b0/-2/0,c0/-1/0,d0/-1/1,a1/0/0,b1/0/1,c1/1/0,d1/1/1,a2/2/0,b2/2/1,c2/3/0}
        \drawedges{b1/d0,d0/b0,b0/c0,c0/a1,a1/c1,c1/b1,b1/d1,d1/a1,c1/a2,a2/c2,c2/b2,b2/d1}
        \draw (b0) to[bend right=28.2] (c2);
        \nodesvertex{c0,d0,b0,a1,b1,c1,d1,a2,b2,c2}
    \end{tikzpicture}
    \caption{The sporadic bipartite graph $\hjone$ of girth $4$ without eigenvalues in $(-1,1)$.} \label{fig:bipartite-girth-4-sporadic}
\end{figure}

\begin{lemma} \label{lem:key}
    For every rooted graph $G_R$, the modified adjacency matrix of $G_R$ has smallest eigenvalue at least $-2$ if and only if, for every $n \in \N$, the graph $G_R + K_n$ has smallest eigenvalue at least $-2$, where the graph $G_R + K_n$ is obtained from $G$ by attaching a copy of the clique $K_n$ of order $n$ to each root in $R$.
\end{lemma}

\begin{proof}
    Let $M$ be the modified adjacency matrix of $G_R$, and let $M_n$ be the matrix obtained from $M$ by increasing the diagonal entries corresponding to the roots in $R$ from $-1$ to $-1 + 1/(n+1)$. We write $\lambda_1(\cdot)$ for the smallest eigenvalue of a matrix.

    We claim that $\lambda_1(M) \ge -2$ if and only if $\lambda_1(M_n) \ge -2$ for all $n \in \N$. Clearly, since $M_n \succeq M$, we know that $\lambda_1(M) \ge -2$ implies that $\lambda_1(M_n) \ge -2$ for all $n \in \N$. Assume that $\lambda_1(M) < -2$, that is, there exists $x \in \mathbb{R}^{V(G)}$ such that $x^\intercal M x < -2x^\intercal x$. By a standard continuity argument, $x^\intercal M_n x < -2x^\intercal x$ for some $n \in \N$, which implies that $\lambda_1(M_n) < -2$ for some $n \in \N$.

    It follows immediately once we show that $\lambda_1(M_n) \ge -2$ if and only if the graph $G_R + K_n$ has smallest eigenvalue at least $-2$. We partition the matrix $A_{G_R + K_n} + 2I$ into the following blocks:
    \[
        \begin{pmatrix}
            A_{G-R} + 2I & B & \\
            B^\intercal & A_{G[R]} + 2I & C \\
            & C^\intercal & D + 2I
        \end{pmatrix},
    \]
    where $D$ is the adjacency matrix of the cliques attached to $R$. Since the smallest eigenvalue of $D$ is $-1$, the block $D + 2I$ is positive definite. Therefore, the above block matrix is positive semidefinite if and only if the Schur complement
    \[
        \begin{pmatrix}
            A_{G-R} + 2I & B\\
            B^\intercal & A_{G[R]} + 2I - C(D+2I)^{-1}C^\intercal
        \end{pmatrix}
    \]
    of $D + 2I$ is positive semidefinite. Let $r \in R$ be an arbitrary vertex of $G[R]$. Since the only nonzero entry of $C$ on the row indexed by $r$ is its $(r, v)$ entry for $v$ in the clique $K_n$ attached to $r$, the matrix $C(D+2I)^{-1}C^\intercal$ is a diagonal matrix, and its $(r,r)$ entry simplifies to $\bm{1}^\intercal (A_{K_n}+2I)^{-1} \bm{1}$, where $\bm{1}$ is the $n$-dimensional all-ones column vector. Since $(A_{K_n} + 2I)^{-1} = I - \bm{1}\bm{1}^\intercal/(n+1)$, we have $\bm{1}^\intercal (A_{K_n}+2I)^{-1} \bm{1} = n - n^2/(n+1) = n/(n+1)$. Therefore
    \[
        A_{G[R]} + 2I - C(D + 2I)^{-1}C^\intercal = A_{G[R]} + 2I - \frac{n}{n+1}I = A_{G[R]} + \left(1+\frac{1}{n+1}\right)I,
    \]
    and so the Schur complement of $D+2I$ is equal to $M_n + 2I$.
\end{proof}

We now recall the relevant notions and two classical results.

\begin{definition}[Graph with petals and generalized line graph]
    A \emph{graph with petals} $F$ is a multigraph obtained from a simple graph by attaching pendant double edges, called \emph{petals}. A \emph{generalized line graph} $L(F)$ is the line graph of a graph with petals $F$ in which two vertices of $L(F)$ are adjacent exactly when the corresponding edges of $F$ share a single vertex.
\end{definition}

\begin{theorem}[Theorem 2.1 of Hoffman \cite{H77}] \label{thm:line-graph-at-least-minus-2}
    If $G$ is a generalized line graph, then the smallest eigenvalue of $G$ is at least $-2$. \qed
\end{theorem}

\begin{theorem}[Theorem 4.2, 4.3 and 4.10 of Cameron et al.~\cite{CGSS76}] \label{thm:cgss}
    For every connected graph $G$ on more than $36$ vertices, if the smallest eigenvalue of $G$ is at least $-2$, then $G$ is a generalized line graph.
    \qed
\end{theorem}

As a consequence, we obtain the following structural characterization of rooted distance-two induced subgraphs $G_R$ of a subcubic graph $H$ with minimum degree $2$.

\begin{corollary} \label{lem:line-graph}
    For every connected rooted graph $G_R$ with $R \neq \varnothing$, the modified adjacency matrix of $G_R$ has smallest eigenvalue at least $-2$ if and only if there exists a connected graph with petals $F$ such that $G = L(F)$ and each root in $R$, viewed as a vertex of $G$, represents a pendant edge of $F$.
\end{corollary}

\begin{proof}
    In view of \cref{lem:key}, it suffices to show that the graph $G_R + K_n$ has smallest eigenvalue at least $-2$ for every $n \in \N$ if and only if $G = L(F)$ for some connected graph with petals $F$, and each root in $R$ represents a pendant edge of $F$. We write $\lambda_1(\cdot)$ for the smallest eigenvalue of a graph.

    Assume that $\lambda_1(G_R + K_n) \ge -2$ for every $n \in \N$. In particular, $\lambda_1(G_R + K_{36}) \ge -2$. By \cref{thm:cgss}, the graph $G_R + K_{36}$, which has more than $36$ vertices, must be the line graph of some graph with petals, say $\tilde{F}$.

    Fix a root $r \in R$, and let $e(r)$ be the edge of $\tilde{F}$ corresponding to $r$. The $36$ vertices of the clique attached to $G_R$ at $r$ represent $36$ edges of $\tilde{F}$ that form a star, which we denote by $S_{36}(r)$. By construction, the only edge of $\tilde{F}$ incident with the star $S_{36}(r)$ is $e(r)$. Hence, upon deleting the star $S_{36}(r)$, the edge $e(r)$ becomes a pendant edge.

    Let $F$ be the graph with petals obtained from $\tilde{F}$ by deleting each of the stars $S_{36}(r)$ for $r \in R$, together with any isolated vertices that arise. One verifies directly that $G = L(F)$, and, for every $r \in R$, the edge $e(r)$ is a pendant edge of $F$. Since $G$ is connected and $F$ has no isolated vertices, $F$ is connected.

    Suppose conversely that $G = L(F)$ for some connected graph with petals $F$, and each $r \in R$ represents a pendant edge $e(r)$ of $F$. Fix an arbitrary $n \in \N$, and construct a multigraph $\tilde{F}$ from $F$ by attaching a copy of the star of size $n$ to each pendant edge $e(r)$ for $r \in R$. One checks that $G_R + K_n$ is then the line graph of $\tilde{F}$. By \cref{thm:line-graph-at-least-minus-2}, it follows that $\lambda_1(G_R + K_n) \ge -2$.
\end{proof}

The preceding spectral preparation similarly allows us to focus on the following purely structural classification of bipartite graphs of girth at least $6$.

\begin{theorem} \label{thm:bipartite-girth-6}
    Let $H$ be a connected subcubic bipartite graph of minimum degree $2$ and girth at least $6$ such that neither \cref{fig:g4-f} nor \cref{fig:g4-g} is a rooted distance-two induced subgraph of $H$. If, for every rooted distance-two component $G_R$ of $H$ with $R \neq \varnothing$, there exists a connected graph with petals $F$ such that $G = L(F)$ and each root in $R$, viewed as a vertex of $G$, represents a pendant edge of $F$, then $H$ is isomorphic to one of the four graphs in \cref{fig:bipartite-girth-6}.
\end{theorem}

\begin{figure}
    \centering
    \begin{subfigure}{3.5cm}
        \centering
        \begin{tikzpicture}[baseline=(A.base)]
            \coordinate (A) at (0,0);
            \coordinate (B) at ($(A)+(30:1)$);
            \coordinate (C) at ($(A)+(-30:1)$);
            \coordinate (D) at ($(C)+(30:2)$);
            \coordinate (E) at ($(B)+(-30:2)$);
            \coordinate (F) at ($(E)+(30:1)$);
            \draw (A)--(B)--(E)--(F)--(D)--(C)--cycle;
            \nodesvertex{A,B,C,D,E,F}
        \end{tikzpicture}
        \caption{} \label{fig:g6-a}
    \end{subfigure}
    \begin{subfigure}{3.5cm}
        \centering
        \begin{tikzpicture}[baseline=(A.base)]
            \coordinate (A) at (0,0);
            \coordinate (B) at ($(A)+(30:1)$);
            \coordinate (C) at ($(A)+(-30:1)$);
            \coordinate (D) at ($(C)+(30:2)$);
            \coordinate (E) at ($(B)+(-30:2)$);
            \coordinate (F) at ($(E)+(30:1)$);
            \coordinate (B1) at ($(B)+(60:.866)$);
            \coordinate (D1) at ($(D)+(120:.866)$);
            \draw (B)--(E)--(F)--(D)--(C)--(A)--(B)--(B1)--(D1)--(D);
            \nodesvertex{A,B,C,D,E,F,B1,D1}
        \end{tikzpicture}
        \caption{} \label{fig:g6-b}
    \end{subfigure}
    \begin{subfigure}{4cm}
        \centering
        \begin{tikzpicture}[baseline=(B.base)]
            \coordinate (A) at (0,0);
            \coordinate (B) at (0,1);
            \coordinate (C) at ($(A)+(150:1)$);
            \coordinate (D) at ($(C)+(90:2)$);
            \coordinate (E) at ($(B)+(150:2)$);
            \coordinate (F) at ($(E)+(90:1)$);
            \coordinate (G) at ($(E)+(-120:.866)$);
            \coordinate (H) at ($(C)+(180:.866)$);
            \coordinate (I) at ($(A)+(30:1)$);
            \coordinate (J) at ($(I)+(90:2)$);
            \coordinate (K) at ($(B)+(30:2)$);
            \coordinate (L) at ($(K)+(90:1)$);
            \coordinate (M) at ($(K)+(-60:.866)$);
            \coordinate (N) at ($(I)+(0:.866)$);
            \draw (E)--(F)--(D)--(C)--(A)--(B)--(E)--(G)--(H)--(C);
            \draw (A)--(I)--(J)--(L)--(K)--(B);
            \draw (I)--(N)--(M)--(K);
            \nodesvertex{A,B,C,D,E,F,G,H,I,J,K,L,M,N}
        \end{tikzpicture}
        \caption{} \label{fig:g6-c}
    \end{subfigure}
    \begin{subfigure}{4cm}
        \centering
        \begin{tikzpicture}[baseline=(B.base)]
            \coordinate (A) at (0,0);
            \coordinate (B) at (0,1);
            \coordinate (C) at ($(A)+(150:1)$);
            \coordinate (D) at ($(C)+(90:2)$);
            \coordinate (E) at ($(B)+(150:2)$);
            \coordinate (F) at ($(E)+(90:1)$);
            \coordinate (G) at ($(E)+(-120:.866)$);
            \coordinate (H) at ($(C)+(180:.866)$);
            \coordinate (I) at ($(A)+(30:1)$);
            \coordinate (J) at ($(I)+(90:2)$);
            \coordinate (K) at ($(B)+(30:2)$);
            \coordinate (L) at ($(K)+(90:1)$);
            \coordinate (M) at ($(K)+(-60:.866)$);
            \coordinate (N) at ($(I)+(0:.866)$);
            \coordinate (O) at ($(D)+(-30:1)$);
            \coordinate (R) at ($(F)+(30:2)$);
            \coordinate (P) at ($(O)!1/3!(R)$);
            \coordinate (Q) at ($(R)!1/3!(O)$);
            \draw (E)--(F)--(D)--(C)--(A)--(B)--(E)--(G)--(H)--(C);
            \draw (A)--(I)--(J)--(L)--(K)--(B);
            \draw (I)--(N)--(M)--(K);
            \draw (F)--(R)--(L);
            \draw (D)--(O)--(J);
            \draw (O)--(P)--(Q)--(R);
            \nodesvertex{A,B,C,D,E,F,G,H,I,J,K,L,M,N,O,P,Q,R}
        \end{tikzpicture}
        \caption{} \label{fig:g6-d}
    \end{subfigure}
    \caption{Sporadic bipartite graphs of girth $6$ without eigenvalues in $(-1,1)$.}
    \label{fig:bipartite-girth-6}
\end{figure}

To reduce the non-bipartite case to the bipartite one, we use the well-known bipartite double construction: if $G$ has eigenvalues $\sset{\lambda_i}$, then its bipartite double has eigenvalues $\sset{\pm \lambda_i}$. This spectral reduction allows us to focus on the following purely structural classification for the non-bipartite case.

\begin{theorem} \label{thm:non-bipartite}
    For every graph $H$, if its bipartite double is isomorphic to $\hjn$ for some $n \ge 1$, $\hjone$, or one of the four graphs in \cref{fig:bipartite-girth-6}, then $H$ is isomorphic to the triangle graph $K_3$ or the graph $\hjn'$ in \cref{fig:non-bipartite} for some $n \ge 1$.
\end{theorem}

Finally, we complete the proof of the classification and derive its consequence on spectral gap sets.

\begin{proof}[Proof of \cref{thm:main} assuming \cref{thm:bipartite-girth-4,thm:bipartite-girth-6,thm:non-bipartite}]
    Let $H$ be a connected subcubic graph that is not cubic and has no eigenvalues in $(-1,1)$. If the minimum degree of $H$ is at most $1$, then \cref{thm:min-deg-1} implies that $H$ is isomorphic to $K_2$. We may therefore assume that the minimum degree of $H$ is $2$.

    We break into the following cases.
    \begin{kase}
        \item Suppose that $H$ is bipartite and has girth $4$. Then \cref{lem:psd,lem:forb-bipartite} imply that none of the rooted multigraphs in \crefrange{fig:g4-a}{fig:g4-f} is a rooted distance-two induced subgraph of $H$, so \cref{thm:bipartite-girth-4} applies.

        \item Suppose that $H$ is bipartite and has girth at least $6$. Then \cref{lem:psd,lem:forb-bipartite} imply that neither \cref{fig:g4-f} nor \cref{fig:g4-g} is a rooted distance-two induced subgraph of $H$. Let $G_R$ be a rooted distance-two component of $H$ with $R \neq \varnothing$. By \cref{lem:psd}, its modified adjacency matrix has smallest eigenvalue at least $-2$, so \cref{lem:line-graph} gives the required representation of $G$ as the line graph of a connected graph with petals. Thus \cref{thm:bipartite-girth-6} applies.

        \item Suppose that $H$ is non-bipartite, and let $H'$ be its bipartite double. The graph $H'$ is connected, subcubic, bipartite, and has minimum degree $2$ and no eigenvalues in $(-1,1)$. Applying the previous two cases to $H'$, we conclude that $H'$ is isomorphic to $\hjn$ for some $n \ge 1$, $\hjone$, or one of the four graphs in \cref{fig:bipartite-girth-6}. Thus \cref{thm:non-bipartite} applies.
    \end{kase}
    We conclude that $H$ is isomorphic to $\hjn$ or $\hjn'$ for some $n \ge 1$, or one of the seven sporadic examples.

    Conversely, we break into the following cases.
    \begin{kase}
        \item Suppose that $H$ is $\hjn$. Both rooted distance-two components of $H$ are as follows.
        \figurehc[1.2]{
            \defxy{a/0/0.5, b/1/0, c/1/1, d/2/0, e/2/1, f/4/0, g/4/1, h/5/0, i/5/1}
            \draw[dashed] (d) to (f);
            \draw[dashed] (e) to (g);
            \draw (b) to[bend right=30] (c) to[bend right=30] (b);
            \draw (d) to[bend right=30] (e) to[bend right=30] (d);
            \draw (f) to[bend right=30] (g) to[bend right=30] (f);
            \drawedges{a/b,a/c,b/d,c/e,f/h,g/i,h/i}
            \node[root-vertex, label=left:$a_0$] at (a) {};
            \node[vertex, label=below:$a_1$] at (b) {};
            \node[vertex, label=above:$a_1'$] at (c) {};
            \node[vertex, label=below:$a_2$] at (d) {};
            \node[vertex, label=above:$a_2'$] at (e) {};
            \node[vertex, label=below:$a_n$] at (f) {};
            \node[vertex, label=above:$a_n'$] at (g) {};
            \node[root-vertex, label=below:$a_{n+1}$] at (h) {};
            \node[root-vertex, label=above:$a_{n+1}'$] at (i) {};
        }
        For $X \subseteq \sset{a_0, a_1, a_1', \dots, a_{n+1}, a_{n+1}'}$, let $J(X) = \bm{1}_X\bm{1}_X^\intercal$, where $\bm{1}_X$ is the indicator vector of $X$. Since the modified adjacency matrix of the above rooted multigraph plus $2I$ is equal to \[
            J(a_0a_1a_1') + J(a_1a_1'a_2a_2') + \dots + J(a_na_n'a_{n+1}a_{n+1}'),
        \]
        it is positive semidefinite, which implies that $H$ has no eigenvalues in $(-1,1)$ by \cref{lem:psd}.
        \item If $H$ is $K_2$, then its eigenvalues are $\pm 1$.
        \item Suppose that $H$ is the graph in \cref{fig:bipartite-girth-4-sporadic} or one of the four graphs in \cref{fig:bipartite-girth-6}. Both rooted distance-two components of $H$ are as follows.
        \figurehc[0.6]
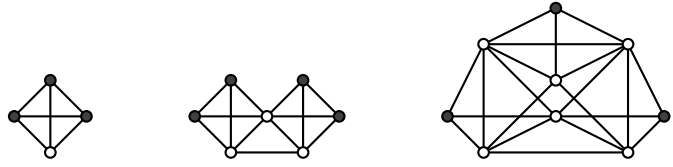
{
            \begin{scope}[shift={(-12,1)}]
                \defxy{a/-0.8660254/0,b/1/1,c/1/-1,d/3/1,e/3/-1}
                \drawedges{a/b,a/c,b/d,b/e,c/d,c/e,a/d,a/e,d/e}
                \nodesrootvertex{d,e}
                \nodesvertex{a,b,c}
            \end{scope}
            \begin{scope}[shift={(-5,-1)}]
                \defxy{a/-1/1,b/0/2,c/1/1}
                \drawedges{a/b,b/c,a/c}
                \nodesrootvertex{a,b,c}
            \end{scope}
            \defxy{a/-1/1,b/0/2,c/1/1,d/0/0}
            \drawedges{a/b,b/c,c/d,d/a,a/c,b/d}
            \nodesrootvertex{a,b,c}
            \nodesvertex{d}
            \begin{scope}[shift={(5,0)}]
                \defxy{a/-1/1,b/0/0,c/1/1,d/0/2,e/2/0,f/3/1,g/2/2}
                \drawedges{a/b,b/c,c/d,d/a,a/c,b/d,c/e,c/f,c/g,e/f,e/g,f/g,b/e}
                \nodesrootvertex{a,d,f,g}
                \nodesvertex{b,c,e}
            \end{scope}
            \begin{scope}[shift={(12,0)}]
                \defxy{a/-1/1,b/0/0,c/2/1,d/0/3,e/4/0,f/5/1,g/4/3,h/2/2,i/2/4}
                \drawedges{a/b,b/c,c/d,d/a,a/c,b/d,c/e,c/f,c/g,e/f,e/g,f/g,b/e,d/h,h/i,d/i,g/h,g/i,b/h,e/h,d/g}
                \nodesrootvertex{a,f,i}
                \nodesvertex{b,c,d,e,g,h}
            \end{scope}
        }
        Since each rooted graph above is the line graph of a graph below in such a way that each root represents a pendant edge, the graph $H$ has no eigenvalues in $(-1,1)$ by \cref{lem:psd,lem:line-graph}.
        \figurehc{
            \begin{scope}[shift={(-5.5,0)}]
                \defpolar{a/0/0,b/0/1,c/60/1,d/120/1,e/180/1}
                \drawedges{a/b,a/c,a/e}
                \draw (a) to[bend right=30] (d);
                \draw (a) to[bend left=30] (d);
            \end{scope}
            \begin{scope}[shift={(-3,0)}]
                \defpolar{a/0/0,b/0/1,c/60/1,e/120/1}
                \drawedges{a/b,a/c,a/e}
            \end{scope}
            \defpolar{a/0/0,b/0/1,c/60/1,d/120/1,e/180/1}
            \drawedges{a/b,a/c,a/d,a/e}
            \begin{scope}[shift={(3,0)}]
                \defpolar{a/0/0,b/0/1,c/0/2,d/180/1,e/60/1,f/120/1}
                \coordinate (g) at ($(e)+(1,0)$);
                \drawedges{a/b,b/c,a/d,f/a,a/e,e/b,b/g}
            \end{scope}
            \begin{scope}[shift={(7,0)}]
                \defpolar{a/0/0,b/0/1,c/0/2,d/180/1,e/60/1,f/120/1}
                \coordinate (g) at ($(e)+(1,0)$);
                \drawedges{a/b,b/c,a/d,f/a,a/e,e/b,b/f,f/e,e/g}
            \end{scope}
        }
        \item Suppose that $H$ is $\hjn'$ or $K_3$. The bipartite double of $H$ is $\mathsf{HJ}_{2n}$ or \cref{fig:g6-a}. Since neither of these graphs has eigenvalues in $(-1,1)$, the same conclusion holds for $H$. \qedhere
    \end{kase}
\end{proof}

To prove the corollary on maximal spectral gap sets, we require the following results.

\begin{lemma}[Theorem 2.1 of Guo and Mohar \cite{GM14}] \label{lem:gmn}
    For $n \in \N^+$, the eigenvalues of $\gmn$ consists of $\pm 1$ each with multiplicity $n$, and $\pm \sqrt{5 + 4\cos (2\pi i/n)}$ for $i \in \sset{0, \dots, n-1}$. \qed
\end{lemma}

\begin{lemma}[Lemma 3.7 of Guo and Royle \cite{GR26}] \label{lem:ksn}
    For $n \in \N^+$, the set of distinct eigenvalues of $\mathsf{GM}_{2n}$ is equal to the set of distinct eigenvalues of $\mathsf{KS}_n$ together with $-3$. \qed
\end{lemma}

\begin{proof}[Proof of \cref{cor:gap}]
    Suppose that $I$ is an open subset of $(-3,-1) \cup (1,3)$. In view of \cref{thm:main}, it suffices to show that for every $n_0 \in \N^+$, there exist $n_1, n_2, n_3, n_4 \ge n_0$ such that each of ${\mathsf{GM}_{n_1}}, {\mathsf{KS}_{n_2}}, {\mathsf{HJ}_{n_3}}, {\mathsf{HJ}_{n_4}'}$ has an eigenvalue in $I$. Fix $n_0 \in \N$. \cref{lem:gmn,lem:ksn} imply that there exist $n_1 \ge n_0$ and $n_2 \ge n_0 + 2$ such that ${\mathsf{GM}_{n_1}}$ has an eigenvalue in $I$, and ${\mathsf{KS}_{n_2}}$ has at least three eigenvalues in $I$. Since $\hjn$ can be obtained from $\mathsf{KS}_{n+2}$ by removing two vertices, and $\hjn'$ can be obtained from $\mathsf{KS}_{n+1}$ by removing one vertex, by the Cauchy interlacing theorem, both ${\mathsf{HJ}_{n_3}}$ and ${\mathsf{HJ}_{n_4}'}$ have eigenvalues in $I$, where $n_3 = n_2-2$ and $n_4 = n_2 - 1$.
\end{proof}

As an independent check of the structural arguments, we developed a Lean 4 formalization of \cref{thm:bipartite-girth-4,thm:bipartite-girth-6,thm:non-bipartite} and the supporting lemmas in \cref{sec:bipartite-girth-4,sec:bipartite-girth-6,sec:non-bipartite} using Codex 0.149 with OpenAI's GPT-5.6 Sol model at high reasoning effort. The formalization led to several improvements in presentation. We reviewed the formal statements and confirmed that the Lean project builds successfully. The Lean source files are included in the ancillary directory \texttt{anc/lean} with the arXiv submission.

The rest of the paper is organized as follows. In \cref{sec:bipartite-girth-4,sec:bipartite-girth-6} we prove the structural classifications in \cref{thm:bipartite-girth-4,thm:bipartite-girth-6} for bipartite graphs of girth $4$ and of girth at least $6$, respectively. In \cref{sec:non-bipartite} we analyze bipartite doubles and prove the structural classification in \cref{thm:non-bipartite}. Finally, in \cref{sec:remark} we discuss further extremal spectral questions, including spectral gap intervals, the spectral gap set $(-2,0)$ and median eigenvalues of subcubic graphs.

\section{\texorpdfstring{Bipartite graphs of girth $4$}{Bipartite graphs of girth 4}} \label{sec:bipartite-girth-4}

In this section, we prove the structural classification in \cref{thm:bipartite-girth-4}.

\begin{proof}[Proof of \cref{thm:bipartite-girth-4}]
    Let $H$ be a connected subcubic bipartite graph of minimum degree $2$ and girth $4$, and suppose that none of the rooted multigraphs in \crefrange{fig:g4-a}{fig:g4-f} is a rooted distance-two induced subgraph of $H$. Let $n \in \N$ be the maximum for which the twisted ladder $\upbowtie_n$ on $n$ rungs is a subgraph of $H$. Denote by $(a,b)$ and $(y,z)$ the two pairs of vertices of the twisted ladder $\upbowtie_n$ as labeled in \cref{fig:twisted-ladder} --- these four vertices are the only ones of degree $2$ in $\upbowtie_n$, and $a$ is at distance $2$ from $b$ while $y$ is at distance $2$ from $z$ in $\upbowtie_n$.

    We deduce from \cref{fig:g4-a} applied to the rooted distance-two subgraph of $H$ induced by $\sset{a,b}$ that both $a$ and $b$ are of degree $3$. Furthermore, we deduce from \cref{fig:g4-b} applied to the rooted distance-two subgraph of $H$ induced by $\sset{a,b}$ that the third neighbors of $a$ and $b$ are distinct.

    \begin{kase}
        \item Suppose that $y$ or $z$ is a third neighbor of $a$ or of $b$. Up to symmetry, we may assume that $y$ is the third neighbor of $a$. Clearly $n \ge 2$. Note that $z$ cannot be the third neighbor of $b$ because otherwise $H$ would become the Guo--Mohar graph $\gmn$, which is cubic. Let $x$ be one of the common neighbors of $y$ and $z$ in $\upbowtie_n$. Depending on whether $n = 2$, the rooted distance-two subgraph of $H$ induced by $\sset{a,b,x}$ is isomorphic to one of \cref{fig:g4-c,fig:g4-d}.
        \figurehc{
            \defxy{a1/0/0, b1/0/1, c1/1/0, d1/1/1, a2/2/0, b2/2/1, c2/3/0, d2/3/1}
            \drawedges{a1/c1, c1/b1, b1/d1, d1/a1, a2/c2, c2/b2, b2/d2, d2/a2, c1/a2, d1/b2}
            \draw (b1) to [bend left=30] (d2);
            \node[vertex, label=left:$b$] at (a1) {};
            \node[vertex, label=left:$a$] at (b1) {};
            \nodesvertex{c1,d1,a2,b2,c2,d2}
            \node[vertex, label=below:$x$] at (a2) {};
            \node[vertex, label=right:$z$] at (c2) {};
            \node[vertex, label=right:$y$] at (d2) {};
            \begin{scope}[shift={(5,0)}]
                \defxy{a1/0/0, b1/0/1, c1/1/0, d1/1/1, a2/2/0, b2/2/1, c2/4/0, d2/4/1, a3/5/0, b3/5/1, c3/6/0, d3/6/1}
                \drawedges{a1/c1, c1/b1, b1/d1, d1/a1, c1/a2, d1/b2, c2/a3, d2/b3, a3/c3, a3/d3, b3/c3, b3/d3}
                \draw (b1) to [bend left=15] (d3);
                \draw[dashed] (a2)--(c2);
                \draw[dashed] (b2)--(d2);
                \node[vertex, label=left:$b$] at (a1) {};
                \node[vertex, label=left:$a$] at (b1) {};
                \nodesvertex{c1,d1,a2,b2,c2,d2,a3,b3,c3,d3}
                \node[vertex, label=right:$z$] at (c3) {};
                \node[vertex, label=right:$y$] at (d3) {};
                \node[vertex, label=below:$x$] at (a3) {};
            \end{scope}
        }
        \item Suppose that neither $y$ nor $z$ is a third neighbor of $a$ or $b$. Let $c$ and $d$ be the (distinct) third neighbors of $a$ and $b$, respectively.
        \begin{kase}
            \item Suppose that the sets of neighbors of $c$ and $d$ outside $\sset{a,b}$ are distinct. Up to symmetry, we may assume that $c$ has another neighbor $e$ that is not adjacent to $d$. Let $f$ and $g$ be the two common neighbors of $a$ and $b$ in $\upbowtie_n$.
            \begin{kase}
                \item If $e$ is adjacent to both $f$ and $g$, then the rooted distance-two subgraph of $H$ induced by $\sset{a,e}$ is isomorphic to \cref{fig:g4-b}.
                \item If $e$ is adjacent to neither $f$ nor $g$, then the rooted distance-two subgraph of $H$ induced by $\sset{a,b,e}$ is isomorphic to \cref{fig:g4-c}.
                \figurehc{
                    \defxy{a1/0/0, b1/0/1, c1/1/0, d1/1/1, a2/2/0, b2/2/1, c2/4/0, d2/4/1, a3/5/0, b3/5/1, c3/6/0, d3/6/1, a0/-1/0, b0/-1/1, c0/-2/1}
                    \drawedges{a1/c1, c1/b1, b1/d1, d1/a1, c1/a2, d1/b2, c2/a3, d2/b3, a3/c3, a3/d3, b3/c3, b3/d3, a0/a1, b0/b1, b0/c0}
                    \draw[dashed] (a2)--(c2);
                    \draw[dashed] (b2)--(d2);
                    \node[vertex, label=above:$e$] at (c0) {};
                    \node[vertex, label=below:$d$] at (a0) {};
                    \node[vertex, label=above:$c$] at (b0) {};
                    \node[vertex, label=below:$b$] at (a1) {};
                    \node[vertex, label=above:$a$] at (b1) {};
                    \node[vertex, label=below:$g$] at (c1) {};
                    \node[vertex, label=above:$f$] at (d1) {};
                    \nodesvertex{a2,b2,c2,d2,a3,b3,c3,d3}
                    \node[vertex, label=below:$z$] at (c3) {};
                    \node[vertex, label=above:$y$] at (d3) {};
                }
                \item If $e$ is adjacent to exactly one of $f$ and $g$, then the rooted distance-two subgraph of $H$ induced by $\sset{a,b,e}$ is isomorphic to \cref{fig:g4-d}.
            \end{kase}
            \item Suppose that the sets of neighbors of $c$ and $d$ outside $\sset{a,b}$ are equal. By the maximality of $n$, we can deduce that $c$ and $d$ are of degree $2$. Let $e$ be the shared second neighbor of $c$ and $d$. Since $H$ is bipartite, we see that $e$ is distinct from $y$ and $z$.
            \figurehc{
                \defxy{a1/0/0, b1/0/1, c1/1/0, d1/1/1, a2/2/0, b2/2/1, c2/4/0, d2/4/1, a3/5/0, b3/5/1, c3/6/0, d3/6/1, a0/-1/0, b0/-1/1, c0/-2/0}
                \drawedges{a1/c1, c1/b1, b1/d1, d1/a1, c1/a2, d1/b2, c2/a3, d2/b3, a3/c3, a3/d3, b3/c3, b3/d3, a0/a1, b0/b1, b0/c0, a0/c0}
                \draw[dashed] (a2)--(c2);
                \draw[dashed] (b2)--(d2);
                \node[vertex, label=below:${\vphantom{b}e}$] at (c0) {};
                \node[vertex, label=below:$d$] at (a0) {};
                \node[vertex, label=above:$c$] at (b0) {};
                \node[vertex, label=below:$b$] at (a1) {};
                \node[vertex, label=above:$a$] at (b1) {};
                \nodesvertex{c1,d1,a2,b2,c2,d2,a3,b3,c3,d3}
                \node[vertex, label=below:$z$] at (c3) {};
                \node[vertex, label=above:$y$] at (d3) {};
            }
            \begin{kase}
                \item Suppose that $y$ or $z$ is a neighbor of $e$. Depending on whether $n = 1$, and whether $\deg_H(z) = 2$, the rooted distance-two subgraph of $H$ induced by $\sset{c,y,z}$ is isomorphic to one of \cref{fig:g4-c,fig:g4-d}.
                \item Suppose that neither $y$ nor $z$ is a neighbor of $e$. Applying the argument symmetrically to $y$ and $z$ in place of $a$ and $b$, we may assume that the third neighbors, denoted $w$ and $x$, of $y$ and $z$ are distinct, both $w$ and $x$ are of degree $2$, and they share a common second neighbor, denoted $v$.
                \figurehc{
                    \defxy{a1/0/0, b1/0/1, c1/1/0, d1/1/1, a2/2/0, b2/2/1, c2/4/0, d2/4/1, a3/5/0, b3/5/1, c3/6/0, d3/6/1, a0/-1/0, b0/-1/1, c0/-2/0, a4/7/0, b4/7/1, a5/8/0}
                    \drawedges{a1/c1, c1/b1, b1/d1, d1/a1, c1/a2, d1/b2, c2/a3, d2/b3, a3/c3, a3/d3, b3/c3, b3/d3, a0/a1, b0/b1, b0/c0, a0/c0, c3/a4, d3/b4, a4/a5, b4/a5}
                    \draw[dashed] (a2)--(c2);
                    \draw[dashed] (b2)--(d2);
                    \node[vertex, label=below:${\vphantom{b}e}$] at (c0) {};
                    \node[vertex, label=below:$d$] at (a0) {};
                    \node[vertex, label=above:$c$] at (b0) {};
                    \node[vertex, label=below:$b$] at (a1) {};
                    \node[vertex, label=above:$a$] at (b1) {};
                    \nodesvertex{c1,d1,a2,b2,c2,d2,a3,b3,c3,d3}
                    \node[vertex, label=below:$z$] at (c3) {};
                    \node[vertex, label=above:$y$] at (d3) {};
                    \node[vertex, label=below:$x$] at (a4) {};
                    \node[vertex, label=above:${\vphantom{y}w}$] at (b4) {};
                    \node[vertex, label=below:$v$] at (a5) {};
                }
                \begin{kase}
                    \item If $\deg_H(e) = \deg_H(v) = 2$, then $H$ is isomorphic to $\hjn$.
                    \item If $ev \in E(H)$ and $n = 1$, then $H$ is isomorphic to $\hjone$.
                    \item Suppose that $\deg_H(e) = 3$ or $\deg_H(v) = 3$, $ev \notin E(H)$, and $n = 1$. Up to symmetry, we may assume that $\deg_H(e) = 3$. Let $f$ be the third neighbor of $e$. Depending on whether $v$ and $f$ share a neighbor, the rooted distance-two subgraph of $H$ induced by $\sset{c,f,v,y}$ is isomorphic to \cref{fig:g4-e} or \cref{fig:g4-f}.
                    \item Suppose that $\deg_H(e) = 3$ or $\deg_H(v) = 3$, and $n \ge 2$. Up to symmetry, we may assume that $\deg_H(e) = 3$. Let $f$ be the third neighbor of $e$. Consider the shortest path of length $2n-1$ in $\upbowtie_n$ between $a$ and $y$. Let $g$ and $h$ be the vertices on that path that are respectively at distance $1$ and $3$ from $a$.
                    \figurehc{
                        \defxy{a1/0/0, b1/0/1, c1/1/0, d1/1/1, a2/2/0, b2/2/1, c2/6/0, d2/6/1, a3/7/0, b3/7/1, c3/8/0, d3/8/1, a0/-1/0, b0/-1/1, c0/-2/0, a4/9/0, b4/9/1, a5/10/0, z0/-3/0, x3/3/0, y3/3/1, x4/4/0, y4/4/1}
                        \drawedges{a1/c1, c1/b1, b1/d1, d1/a1, c1/a2, d1/b2, c2/a3, d2/b3, a3/c3, a3/d3, b3/c3, b3/d3, a0/a1, b0/b1, b0/c0, a0/c0, c3/a4, d3/b4, a4/a5, b4/a5, c0/z0, a2/x3, b2/y3, a2/y3, b2/x3, x3/x4, y3/y4}
                        \draw[dashed] (x4)--(c2);
                        \draw[dashed] (y4)--(d2);
                        \node[vertex, label=below:$f$] at (z0) {};
                        \node[vertex, label=below:$\vphantom{b}e$] at (c0) {};
                        \node[vertex, label=below:$d$] at (a0) {};
                        \node[vertex, label=above:$\vphantom{g}c$] at (b0) {};
                        \node[vertex, label=below:$b$] at (a1) {};
                        \node[vertex, label=above:$\vphantom{g}a$] at (b1) {};
                        \node[vertex, label=above:$g$] at (d1) {};
                        \node[vertex, label=above:$h$] at (y3) {};
                        \nodesvertex{c1,a2,b2,c2,d2,a3,b3,c3,d3,x3,x4,y4}
                        \node[vertex, label=below:$z$] at (c3) {};
                        \node[vertex, label=above:$y$] at (d3) {};
                        \node[vertex, label=below:$x$] at (a4) {};
                        \node[vertex, label=above:${\vphantom{y}w}$] at (b4) {};
                        \node[vertex, label=below:$v$] at (a5) {};
                    }
                    \figurehc{
                        \defxy{b0/-2/0,c0/-1/0,d0/-1/1,a1/0/0,b1/0/1,c1/1/0,d1/1/1,a2/2/0,b2/2/1,c2/3/0,d2/3/1,a3/4/0,b3/4/1,c3/5/0}
                        \drawedges{b1/d0,d0/b0,b0/c0,c0/a1,a1/c1,c1/b1,b1/d1,d1/a1,a2/c2,c2/b2,b2/d2,d2/a2,c1/a2,d1/b2,d2/b3,b3/c3,c3/a3,a3/c2}
                        \draw (b0) to[bend right=30] (c3);
                        \nodesvertex{a1,b1,c1,d1,a2,b2,c2,d2,a3,b3,c3,b0,c0,d0}
                        \node[vertex, label=left:$e$] at (b0) {};
                        \node[vertex, label=above:$d$] at (c0) {};
                        \node[vertex, label=above:$\vphantom{g}c$] at (d0) {};
                        \node[vertex, label=above:$b$] at (a1) {};
                        \node[vertex, label=above:$\vphantom{g}a$] at (b1) {};
                        \node[vertex, label=above:$g$] at (d1) {};
                        \node[vertex, label=above:$\vphantom{g}h$] at (d2) {};
                        \node[vertex, label=above:$\vphantom{g}w$] at (b3) {};
                        \node[vertex, label=right:$f$] at (c3) {};
                        \node[vertex, label=above:$x$] at (a3) {};
                        \node[vertex, label=above:$z$] at (c2) {};
                    }
                    Depending on whether $f$ and $h$ share a neighbor, the rooted distance-two subgraph of $H$ induced by $\sset{c,f,g,h}$ is isomorphic to \cref{fig:g4-e} or \cref{fig:g4-f}. \qedhere
                \end{kase}
            \end{kase}
        \end{kase}
    \end{kase}
\end{proof}

\section{\texorpdfstring{Bipartite graphs of girth at least $6$}{Bipartite graphs of girth at least 6}} \label{sec:bipartite-girth-6}

Observe that if $H$ is a connected subcubic bipartite graph, then its distance-two graph has exactly two connected components. Our strategy is to encode one component in the other, allowing us to focus primarily on a single component. To this end, we introduce the following notions.

\begin{definition}[Valid decomposition, closure, intersection graph, and incidence graph]
    Let $G_R$ be a rooted graph. A \emph{valid decomposition} $\D$ of $G_R$ is a triangle-edge decomposition of $G$ such that every vertex of $G$ belongs to exactly $3$ parts of $\D$, except that every root in $R$ belongs to exactly $2$ parts. The \emph{closure} of an edge $ab$ of $G$, denoted by $\ol{ab}$, is the unique part of $\D$ containing $ab$. The \emph{intersection graph} of $\D$, denoted $I(\D)$, is the graph with vertex set $\D$ in which two vertices are adjacent if and only if their intersection is nonempty. The \emph{incidence graph} of $\D$ is the bipartite graph with parts $\cup \D$ and $\D$, where $a \in \cup \D$ and $\alpha \in \D$ are adjacent if and only if $a \in \alpha$.
\end{definition}

Given a valid decomposition $\D$ of a rooted graph, we denote by $\D_{\sharp 2}$ the set of edge parts and by $\D_{\sharp 3}$ the set of triangle parts in $\D$. We now encode one rooted distance-two component into the other via a valid decomposition.

\begin{lemma} \label{lem:valid-decomposition}
    Let $H$ be a connected subcubic bipartite graph with minimum degree $2$, and let $G_R$ and $G'_S$ be the rooted distance-two components of $H$. If the girth of $H$ is at least $6$, then there exists a valid decomposition $\D$ of $G_R$ such that the rooted graphs $G'_S$ and $\iddtwo$ are isomorphic, and $H$ and the incidence graph of $\D$ are isomorphic.
\end{lemma}

\begin{proof}
    Let $A$ and $B$ be the vertex sets of $G_R$ and $G'_S$ respectively. Since the girth of $H$ is at least $6$, both $G$ and $G'$ are simple graphs. Define $\D = \dset{N_H(b)}{b \in B}$, where $N_H(b)$ denotes the neighborhood of $b$ in $H$. One readily check that $\D$ is a valid decomposition of $G_R$, that the mapping $b \mapsto N_H(b)$ defines an isomorphism from $G'$ to $I(\D)$ and a bijection from $S$ to $\D_{\sharp 2}$, and that $H$ is isomorphic to the incidence graph of $\D$.
\end{proof}

\begin{example}
    Consider the following connected subcubic bipartite graph $H$ of girth $6$, where the vertices in one part are labeled by $a,b,c,d$.
    \figurehc{
        \coordinate (A) at (0,0);
        \coordinate (B) at ($(A)+(30:1)$);
        \coordinate (C) at ($(A)+(-30:1)$);
        \coordinate (D) at ($(C)+(30:2)$);
        \coordinate (E) at ($(B)+(-30:2)$);
        \coordinate (F) at ($(E)+(30:1)$);
        \coordinate (B1) at ($(B)+(60:.866)$);
        \coordinate (D1) at ($(D)+(120:.866)$);
        \draw (B)--(E)--(F)--(D)--(C)--(A)--(B)--(B1)--(D1)--(D);
        \nodesvertexbadgehollow{A/a,D/d,E/b,B1/c}
    }
    The rooted distance-two component $G_R$ induced by $\sset{a,b,c,d}$ is the complete graph $K_4$ with $R = \sset{a,b,c}$. Note that the vertices in another part of $H$ have neighbors $abc, ad, bd, cd$, which form a valid decomposition of $G_R$.
\end{example}

Let $H$ be a connected subcubic bipartite graph of minimum degree $2$ and girth at least $6$, and let $G_R$ and $G'_S$ be its rooted distance-two components, where $R \neq \varnothing$. Suppose that $G$ is the line graph of a connected graph with petals, say $F$, in such a way that each vertex in $R$ represents a pendant edge of $F$. By \cref{lem:valid-decomposition}, there exists a valid decomposition $\D$ of $G_R$ such that $G_S'$ and $\iddtwo$ are isomorphic.

Our analysis therefore reduces to understanding the structure of $F$ together with the decomposition $\D$. We first show that $F$ has no petal, and hence is a simple graph.

\begin{lemma} \label{lem:no-petal}
    Let $F$ be a connected graph with petals and let $R$ be a nonempty set of pendant edges of $F$. If the rooted graph $L(F)_R$ admits a valid decomposition $\D$, and neither \cref{fig:g4-f} nor \cref{fig:g4-g} is an induced rooted subgraph of the rooted graph $\iddtwo$, then $F$ has no petal.
\end{lemma}

The proof proceeds by a case analysis around a petal in $F$. To restrict the possible configurations, we establish the following auxiliary result.

\begin{proposition} \label{lem:3-petals}
    Let $G_0$ be the graph obtained from the complete graph $K_n$ of order $n$ by removing a matching of size $m \in \N^+$, and let $\D_0$ be a triangle-edge decomposition of $G_0$. If $n \in \sset{6,7}$, and each vertex of $G_0$ belongs to $3$ parts of $\D_0$, then $m = 3$.
\end{proposition}

\begin{proof}
    Let $\alpha$ and $\beta$ denote the number of triangle parts and edge parts of $\D_0$ respectively. Counting the vertex-part incidences in two ways gives $3\alpha + 2\beta = 3n$, and hence $\beta$ is divisible by $3$. Counting edges of $G_0$ in two ways yields $\binom{n}{2} - m = 3\alpha + \beta$, which is therefore divisible by $3$. Since $\binom{n}{2}$ is divisible by $3$ and $1 \le m \le n / 2$, we conclude that $m = 3$.
\end{proof}

\begin{proof}[Proof of \cref{lem:no-petal}]
    Set $G = L(F)$. Let $\D$ be a valid decomposition of $G_R$, and set $G' = I(\D)$ and $S = \D_{\sharp 2}$. Assume for the sake of contradiction that $F$ has a petal consisting of two edges $a$ and $b$ with endpoints $u$ and $v$ such that $\deg_F(u) = 2$.

    Since $a \notin R$, it belongs to $3$ parts of $\D$, and therefore $3 \le \deg_G(a) \le 6$. Since $ab \notin E(G)$ and every edge of $F$ incident with $v$ other than $a$ and $b$ is adjacent to $a$ in $G$, it follows that $5 \le \deg_F(v) \le 8$.

    \begin{kase}
        \item Suppose that $\deg_F(v) = 5$. Let $c, d, e$ be the other three edges incident with $v$ in $F$. Since $a \notin R$ and $ac, ad, ae$ are the only edges of $G$ incident with $a$, it follows that $ac, ad, ae \in \D$. Similarly, $bc, bd, be \in \D$. The subgraph of $G'_S$ induced by $\sset{ac, bc, bd}$ is isomorphic to \cref{fig:g4-g}.

        \item Suppose that $\deg_F(v) = 6$. Let $c,d,e,f$ be the other four edges incident with $v$, and set $U = \sset{a,b,c,d,e,f}$. We claim that the edges in $U$ pair up to form three petals in $F$.

        Let $m$ be the number of pairs of edges in $U$ that form petals in $F$, let $G_0$ be the subgraph of $G$ induced by $U$, and let $\D_0$ be the restriction of $\D$ to $G_0$. Then $G_0$ is obtained from $K_6$ by removing a matching of size $m$, and each vertex in $U$ belongs to $3$ parts of $\D_0$ since either it is not in $R$ or its degree in $G_0$  is at least $5$. By \cref{lem:3-petals}, we have $m = 3$.

        Thus, up to symmetry, $\sset{c,d}$ and $\sset{e,f}$ are the other two petals in $F$.
        \figurehc[1.4]{
            \defpolar{a/0/0,b/0/1,c/90/1,d/180/1}
            \draw (a) to [bend right=30] node[midway, fill=white, inner sep=2pt] {$a$} (b);
            \draw (b) to [bend right=30] node[midway, fill=white, inner sep=2pt] {$b$} (a);
            \draw (a) to [bend right=30] node[midway, fill=white, inner sep=2pt] {$c$} (c);
            \draw (c) to [bend right=30] node[midway, fill=white, inner sep=2pt] {$d$} (a);
            \draw (a) to [bend right=30] node[midway, fill=white, inner sep=2pt] {$e$} (d);
            \draw (d) to [bend right=30] node[midway, fill=white, inner sep=2pt] {$f$} (a);
        }
        Since $a \notin R$ and $\deg_G(a) = 4$, the vertex $a$ belongs to exactly one triangle part and two edge parts of $\D$; the same holds for each of $a, b, c, d, e, f$. Without loss of generality, assume that $ace, bdf, ad, af, bc, be, cf, de \in \D$. Then the subgraph of $G'_S$ induced by $\sset{ad, af, cf}$ is isomorphic to \cref{fig:g4-g}.

        \item Suppose that $\deg_F(v) \in \sset{7,8}$. Let $c,d,e,f,g$ be five other edges of $F$ incident with $v$, and set $U = \sset{a,b,c,d,e,f,g}$. As above, let $m$ be the number of petals among pairs in $U$, let $G_0$ be the subgraph of $G$ induced by $U$, and let $\D_0$ be the restriction of $\D$ to $G_0$. Then $G_0$ is obtained from $K_7$ by removing a matching of size $m$, and each vertex in $U$ belongs to $3$ parts in $\D_0$ since its degree in $G_0$ is at least $5$. By \cref{lem:3-petals}, we have $m = 3$.

        Thus, up to symmetry, assume that $\sset{c,d}$ and $\sset{e,f}$ are the other two petals in $F$.
        \figurehc[1.4]{
            \defpolar{a/0/0,b/0/1,c/90/1,d/180/1,e/-90/1}
            \draw (a) to [bend right=30] node[midway, fill=white, inner sep=2pt] {$a$} (b);
            \draw (b) to [bend right=30] node[midway, fill=white, inner sep=2pt] {$b$} (a);
            \draw (a) to [bend right=30] node[midway, fill=white, inner sep=2pt] {$c$} (c);
            \draw (c) to [bend right=30] node[midway, fill=white, inner sep=2pt] {$d$} (a);
            \draw (a) to [bend right=30] node[midway, fill=white, inner sep=2pt] {$e$} (d);
            \draw (d) to [bend right=30] node[midway, fill=white, inner sep=2pt] {$f$} (a);
            \drawlabeledges{a/e/g}
        }
        Since $\deg_G(g) \ge 6$ and $g$ belongs to at most $3$ parts of $\D$, the only edges of $G$ incident with $g$ are $ag, bg, cg, dg, eg, fg$, and $g$ belongs to $3$ triangle parts of $\D$. Since $ab, cd, ef \notin E(G)$, we may assume that these triangle parts are $adg, beg, cfg \in \D$.

        \begin{kase}
            \item Suppose that $\deg_F(v) = 7$. Since $a \notin R$, $\deg_G(a) = 5$, $ab, ef \notin E(G)$, and $adg, cfg \in \D$, it follows that $ace, af \in \D$. Similarly, $bdf, bc \in \D$. Then the subgraph of $G'_S$ induced by $\sset{beg, ace, af, bdf}$ is isomorphic to \cref{fig:g4-f}.

            \item Suppose that $\deg_F(v) = 8$, and let $h$ be the remaining edge of $F$ incident with $v$. Since $gh \notin E(G)$, the edges $g$ and $h$ share both endpoints and form a fourth petal. Since $F$ is connected, $F$ consists of only the eight edges $a,b,c,d,e,f,g,h$, none of which is a pendant edge. \qedhere
        \end{kase}
    \end{kase}
\end{proof}

The next lemma allows us to conclude that the other rooted distance-two component $G_S'$ has at least one root, that is, $S \neq \varnothing$. This, in turn, allows us to apply our arguments symmetrically to the two components $G_R$ and $G'_S$.

\begin{lemma} \label{lem:has-edge-parts}
    Let $F$ be a graph and let $R$ be a nonempty subset of pendant edges of $F$. If the rooted graph $L(F)_R$ admits a valid decomposition $\D$, then $\D$ contains at least one edge part.
\end{lemma}

\begin{proof}
    Set $G = L(F)$. Assume for the sake of contradiction that a valid decomposition $\D$ of $G_R$ consists only of triangle parts. Choose $a \in R$, and let $v$ be the endpoint of $a$ whose other endpoint is a leaf vertex. Since $a \in R$, it belongs to exactly $2$ triangle parts, say $abc$ and $ade$, of $\D$. Consequently, $ab,ac,ad,ae$ are the only edges of $G$ incident with $a$, and hence $a,b,c,d,e$ are the only edges of $F$ incident with $v$.

    Now $bd, be \in E(G)$. Let $f$ and $g$ be vertices of $G$ such that $\ol{bd} = bdf$ and $\ol{be} = beg$. Since $ba,bc,bd,be,bf,bg \in E(G)$ and $b$ belongs to at most $3$ parts of $\D$, it follows that $\deg_G(b) = 6$, and hence $ba, bc, bd, be, bf, bg$ are the precisely the edges of $G$ incident with $b$. Let $u_b, u_d, u_e$ denote the endpoints of $b,d,e$, respectively, other than $v$. Since neither $f$ nor $g$ can be incident with $v$ in $F$, we must have $f = u_bu_d$ and $g = u_bu_e$.

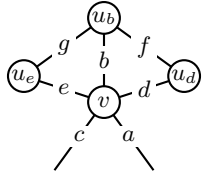
    \figurehc[1.4]{
        \defpolar{v/0/0,u/-54/1,a/18/1,b/90/1,c/162/1,d/234/1}
        \drawlabeledges{u/v/a,v/a/d,v/b/b,v/c/e,v/d/c,a/b/f,b/c/g}
        \nodeslabelvertex{v/v,a/u_d,b/u_b,c/u_e}
    }

    Since $b$ already belongs to three parts $abc,bdf, beg$ of $\D$, the edges $b, f, g$ are the only edges of $F$ incident with $u_b$, and hence $\deg_F(u_b) = 3$. A symmetric argument shows that $\deg_F(u_d) = 3$. Since $\deg_F(u_b) = \deg_F(u_d) = 3$, it follows that $\deg_G(f) = 4$, and therefore $f$ belongs to exactly $2$ triangle parts of $\D$, which contradicts the fact that $f \notin R$.
\end{proof}

Together, \cref{lem:no-petal} and \cref{lem:has-edge-parts} imply that both $G_R$ and $G'_S$ are line graphs, that is, $G_R = L(F)$ and $G'_S = L(F')$ for some connected simple graphs $F$ and $F'$. We now recall two classical results on line graphs.

\begin{theorem}[Whitney~\cite{W32}] \label{thm:whitney}
    If the line graphs of two connected graphs $G$ and $H$ are isomorphic, then $G$ and $H$ are isomorphic, except in the case where the common line graph is isomorphic to the triangle graph $K_3$.
\end{theorem}

\begin{theorem}[Beineke~\cite{B70}] \label{thm:beineke}
    A graph $G$ is a line graph if and only if none of the nine graphs in \cref{fig:beineke} occurs as an induced subgraph of $G$.
\end{theorem}

\begin{figure}
    \centering
    \beineke{a}{
        \defpolar{0/0/0,1/0/1,2/60/1,3/-60/1}
        \drawedges{0/1,0/2,0/3}
        \nodesvertex{0,1,2,3}
    }
    \beineke{b}{
        \defpolar{0/0/0,1/0/1,2/120/1,3/-120/1,4/180/1}
        \drawedges{1/2,0/2,0/3,1/3,0/4,2/4,3/4}
        \nodesvertex{0,1,2,3,4}
    }
    \beineke{c}{
        \defpolar{0/0/0,1/0/1,2/120/1,3/240/1,4/180/1}
        \drawedges{0/1,1/2,0/2,0/3,0/4,1/3,2/3,2/4,3/4}
        \nodesvertex{0,1,2,3,4}
    }
    \beineke{d}{
        \defpolar{0/0/0,1/60/1,2/120/1,3/180/1,4/240/1,5/300/1}
        \drawedges{0/2,1/2,0/3,0/4,4/5,2/3,3/4}
        \nodesvertex{0,1,2,3,4,5}
    }
    \beineke{e}{
        \defpolar{0/0/0,1/0/1,2/120/1,3/240/1,4/180/1,5/180/2}
        \drawedges{0/1,0/2,0/3,2/3,2/4,3/4,2/5,3/5,4/5}
        \nodesvertex{0,1,2,3,4,5}
    }\\\vspace{1em}%
    \beineke{f}{
        \defpolar{0/0/0,1/0/1,2/120/1,3/240/1,4/180/1,5/180/2}
        \drawedges{0/1,1/2,0/2,0/3,1/3,2/3,2/4,3/4,2/5,3/5,4/5}
        \nodesvertex{0,1,2,3,4,5}
    }
    \beineke{g}{
        \defpolar{0/0/0,1/60/1,2/120/1,3/180/1,4/240/1,5/300/1}
        \drawedges{0/2,1/2,0/3,0/4,4/5,2/3,3/4,1/5}
        \nodesvertex{0,1,2,3,4,5}
    }
    \beineke{h}{
        \defpolar{0/0/0,1/60/1,2/120/1,3/180/1,4/240/1}
        \coordinate (5) at (210:{sqrt(3)});
        \drawedges{0/2,1/2,0/3,0/4,4/5,2/3,3/4,0/1,3/5}
        \nodesvertex{0,1,2,3,4,5}
    }
    \beineke{i}{
        \defpolar{0/0/0,1/90/1,2/162/1,3/234/1,4/306/1,5/378/1}
        \drawedges{0/1,0/2,0/3,0/4,0/5,1/2,2/3,3/4,4/5,5/1}
        \nodesvertex{0,1,2,3,4,5}
    }
    \caption{The nine minimal forbidden induced subgraphs for line graphs.} \label{fig:beineke}
\end{figure}

The next two lemmas further characterize $F$ as chordal and free of the diamond graph.

\begin{lemma} \label{lem:chordal}
    Let $F$ be a graph and let $R$ be a subset of pendant edges of $F$. If the rooted graph $L(F)_R$ admits a valid decomposition $\D$, and there exists a graph $F'$ such that $I(\D) = L(F')$ and each edge part of $\D$, viewed as a vertex of $I(\D)$, represents a pendant edge of $F'$, then $F$ is chordal.
\end{lemma}

\begin{proof}
    Set $G = L(F)$, and let $\D$ be a valid decomposition of $G_R$. Assume for the sake of contradiction that $F$ contains an induced cycle $v_1 \dots v_\ell$ of length $\ell \ge 4$. Let its edges be $a_i = v_{i-1} v_i$ for $i \in \sset{1,\dots,\ell}$. Then $a_ia_{i+1} \in E(G)$ for $i \in \sset{1,\dots, \ell}$. Throughout the proof, all the indices are taken modulo $\ell$.
    \figurehc[1.8]{
        \defxy{a/0/0,b/1/0,c/2/0,d/3/0,e/4/0}
        \drawlabeledges{a/b/a_1,b/c/a_2,c/d/a_3,d/e/a_4}
        \nodeslabelvertex{a/v_\ell,b/v_1,c/v_2,d/v_3,e/v_4}
        \draw[dashed, shorten <=12pt, shorten >=12pt] (a) to[bend left=30] (e);
    }

    We claim that if $\ol{a_ia_{i+1}}$ and $\ol{a_ja_{j+1}}$ intersect, then $i \in \sset{j-1, j, j+1}$. Indeed, if they intersect but $\sset{a_i,a_{i+1}} \cap \sset{a_j,a_{j+1}} = \varnothing$, then $\ol{a_ia_{i+1}} = a_ia_{i+1}b$ and $\ol{a_ja_{j+1}} = a_ja_{j+1}b$ for some vertex $b$ of $G$. Since $a_ib, a_{i+1}b \in E(G)$, and $v_1 \dots v_\ell$ is an induced cycle of $F$, the edge $b$ must be incident with $v_i$ in $F$; similarly, $b$ is incident with $v_j$ in $F$. Thus $b$ is a chord of the induced cycle, a contradiction.

    Hence the closures $\ol{a_ia_{i+1}}$ form an induced cycle of length $\ell$ in $I(\D)$. By \cref{thm:whitney}, since $\ell \ge 4$, this induced cycle represents a cycle in $F'$. In particular, none of the vertices $\ol{a_ia_{i+1}}$ represents a pendant edge of $F'$, and therefore $\ol{a_ia_{i+1}}$ is a triangle part of $\D$.

    Let $b_i$ be a vertex of $G$ such that $\ol{a_ia_{i+1}} = a_ib_ia_{i+1}$ for $i \in \sset{1,\dots, \ell}$. Since $v_1 \dots v_\ell$ is an induced cycle of $F$, the edge $b_i$ is incident with $v_i$ in $F$. Let $u_i$ be its other endpoint.

    \begin{kase}
        \item Suppose that all $u_i$ coincide. Since the vertex $b_1$ is adjacent in $G$ to $a_1,a_2,b_2,\dots,b_\ell$, and it belongs to at most $3$ parts of $\D$, we have $\ell + 1 \le \deg_G(b_1) \le 6$, and so $\ell \in \sset{4,5}$.
        \begin{kase}
            \item Suppose that $\ell = 4$. Set $U = \sset{b_1,\dots,b_4}$.
            \figurehc[1.4]{
                \defpolar{a/-45/1,b/45/1,c/135/1,d/-135/1,u/0/0}
                \drawlabeledges{a/b/a_1,b/c/a_2,c/d/a_3,d/a/a_4,u/a/b_4,u/b/b_1,u/c/b_2,u/d/b_3}
            }
            We claim that each $b_i$ belongs to a triangle in $G[U]$ that is a part of $\D$. Without loss of generality, we consider the vertex $b_1$ of $G$. Since $a_1b_1a_2 \in \D$, $b_1b_2, b_1b_3, b_1b_4 \in E(G)$, and $b_1$ belongs to at most $3$ parts of $\D$, we know that $b_1$ belongs to a triangle in $G[U]$. This finishes the proof of the claim. However, $G[U]$, which is $K_4$, has no edge-disjoint triangles that covering all vertices, a contradiction.

            \item Suppose that $\ell = 5$. Set $U = \sset{b_1,\dots,b_5}$.
            \figurehc[1.4]{
                \defpolar{a/18/1,b/90/1,c/162/1,d/234/1,e/306/1,u/0/0}
                \drawlabeledges{a/b/a_1,b/c/a_2,c/d/a_3,d/e/a_4,e/a/a_5,u/a/b_5,u/b/b_1,u/c/b_2,u/d/b_3,u/e/b_4}
            }
            We claim that edges of $G[U]$ are partitioned into triangle parts of $\D$. Without loss of generality, we consider the edge $b_1b_2$ of $G[U]$. Since $\deg_G(b_1) \ge 6$, we know that $b_1$ belongs to $3$ triangle parts of $\D$. Since $a_1b_1a_2 \in \D$, the other two triangle parts incident with $b_1$ are completely contained in $G[U]$. This finishes the proof of the claim. However, $G[U]$ has $10$ edges, which is not divisible by $3$, again a contradiction.
        \end{kase}
        \item Suppose that there exists an index $i$ satisfying $u_{i} \neq u_{i+1}$ and $u_{i+1} \neq u_{i+2}$. Without loss of generality, we may assume that $u_1 \neq u_2$ and $u_2 \neq u_3$.
        \figurehc[1.4]{
            \defxy{a/0/0,b/1/0,c/2/0,d/3/0,e/4/0,b1/1/1,c1/2/1,d1/3/1}
            \drawlabeledges{a/b/a_1,b/c/a_2,c/d/a_3,d/e/a_4,b/b1/b_1,c/c1/b_2,d/d1/b_3}
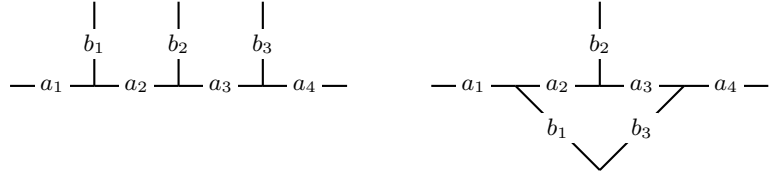
            \begin{scope}[shift={(5,0)}]
                \defxy{a/0/0,b/1/0,c/2/0,d/3/0,e/4/0,b1/2/-1,c1/2/1}
                \drawlabeledges{a/b/a_1,b/c/a_2,c/d/a_3,d/e/a_4,b/b1/b_1,c/c1/b_2,d/b1/b_3}
            \end{scope}
        }
        Since $a_2b_2a_3 \in \D$, and $b_2$ belongs to at least $2$ parts of $\D$, there exists a part $\sigma$ of $\D$ such that $b_2 \in \sigma$, but $a_2, a_3 \notin \sigma$. Since $b_1, a_1, b_3, a_4$ are not adjacent to $b_2$ in $G$, we also have $b_1, a_1, b_3, a_4 \notin \sigma$. Then the subgraph of $I(\D)$ induced by $\sset{a_1b_1a_2, a_2b_2a_3, a_3b_3a_4, \sigma}$ is isomorphic to \cref{fig:beineke-a}, which contradicts \cref{thm:beineke}.
        \figurehc[1.2]{
            \defpolar{0/0/0,1/0/1,2/60/1,3/-60/1}
            \drawedges{0/1,0/2,0/3}
            \nodesvertexbadge{0/a_2b_2a_3,1/\sigma,2/a_1b_1a_2,3/a_3b_3a_4}
        }
        \item Suppose that there exists an index $i$ such that $u_i = u_{i+1}$, $u_{i+2} = u_{i+3}$ but $u_{i+1} \neq u_{i+2}$. Without loss of generality, we may assume that $u_1 = u_2$, $u_3 = u_4$ but $u_2 \neq u_3$.
        \figurehc[1.6]{
            \defxy{a/0/0,b/1/0,c/2/0,d/3/0,e/4/0,f/5/0,b1/1.5/1,c1/3.5/1}
            \drawlabeledges{a/b/a_1,b/c/a_2,c/d/a_3,d/e/a_4,e/f/a_5,b/b1/b_1,c/b1/b_2,d/c1/b_3,e/c1/b_4}
            \nodeslabelvertex{c/v_2}
        }
        Since $a_3 \notin R$ and $a_2b_2a_3, a_3b_3a_4 \in \D$, the vertex $a_3$ belongs to a third part $\ol{a_3c}$ of $\D$ for some $c \notin \sset{a_2, b_2, b_3, a_4}$. Up to symmetry, we may assume $c$ is incident with $v_2$ in $F$. Since $a_2b_2a_3 \in \D$ and $a_2c, b_2c, a_3c \in E(G)$, the three closures $\ol{a_2c}, \ol{b_2c}, \ol{a_3c}$ are three distinct parts of $\D$. Then the subgraph of $I(\D)$ induced by $\sset{\ol{a_2c}, \ol{b_2c}, \ol{a_3c}, a_2b_2a_3, a_3b_3a_4, a_4b_4a_5}$ is isomorphic to \cref{fig:beineke-e}, which contradicts \cref{thm:beineke}.
        \figurehc[1.6]{
            \defpolar{0/0/0,1/0/1,2/120/1,3/240/1,4/180/1,5/180/2}
            \drawedges{0/1,0/2,0/3,2/3,2/4,3/4,2/5,3/5,4/5}
            \nodesvertexbadge{0/a_3b_3a_4,1/a_4b_4a_5,2/\ol{a_3c},3/a_2b_2a_3,4/\ol{b_2c},5/\ol{a_2c}}
        }
        Even when $\ell = 4$ and $a_5 = a_1$, the above argument works exactly the same way. \qedhere
    \end{kase}
\end{proof}

\begin{lemma} \label{lem:diamond-free}
    Let $F$ be a graph and let $R$ be a subset of pendant edges of $F$. If the rooted graph $L(F)_R$ admits a valid decomposition $\D$, and there exists a graph $F'$ such that $I(\D) = L(F')$ and each edge part of $\D$, viewed as a vertex of $I(\D)$, represents a pendant edge of $F'$, then $F$ does not contain the diamond graph as an induced subgraph.
\end{lemma}

\begin{proof}
    Set $G = L(F)$. Assume for the sake of contradiction that $F$ contains an induced diamond, shown below.
    \figurehc[1.4]{
        \defpolar{1/0/0,2/30/1,3/90/1,4/150/1}
        \drawlabeledges{1/2/e,1/3/c,1/4/b,2/3/d,3/4/a}
    }
    Set $U = \sset{a,b,c,d,e}$. Since $\deg_{G[U]}(c) = 4$ and $c$ belongs to at most $3$ parts of $\D$, it belongs to one or two triangles in $G[U]$ that are parts of $\D$.
    \begin{kase}
        \item Suppose that $c$ belongs to exactly one triangle in $G[U]$ that is a part of $\D$.
        \begin{kase}
            \item Suppose that $abc \in \D$ or $cde \in \D$. Up to symmetry, we may assume that $abc \in \D$ and $cde \notin \D$. One can check that $\ol{ad}, \ol{be}, \ol{cd}, \ol{ce}, \ol{de}$ are five distinct parts of $\D$. Then the subgraph of $I(\D)$ induced by $\sset{\ol{ad}, \ol{be}, \ol{ce}, \ol{de}, abc}$ is isomorphic to \cref{fig:beineke-b}, which contradicts \cref{thm:beineke}.
            \figurehc[1.4]{
                \defpolar{0/0/0,1/0/1,2/120/1,3/-120/1,4/180/1}
                \drawedges{1/2,0/2,0/3,1/3,0/4,2/4,3/4}
                \nodesvertexbadge{0/\ol{ce},1/\ol{ad},2/\ol{de},3/abc,4/\ol{be}}
            }
            \item Suppose that $acd \in \D$ or $bce \in \D$. Up to symmetry, we may assume that $acd \in \D$ and $bce \notin \D$. One can check that $\ol{ab}, \ol{bc}, \ol{be}, \ol{ce}, \ol{de}$ are five distinct parts of $\D$. Then the subgraph of $I(\D)$ induced by $\sset{\ol{ab}, \ol{be}, \ol{ce}, \ol{de}, acd}$ is isomorphic to \cref{fig:beineke-b}, which contradicts \cref{thm:beineke}.
            \figurehc[1.4]{
                \defpolar{0/0/0,1/0/1,2/120/1,3/-120/1,4/180/1}
                \drawedges{1/2,0/2,0/3,1/3,0/4,2/4,3/4}
                \nodesvertexbadge{0/\ol{de},1/\ol{ab},2/\ol{be},3/acd,4/\ol{ce}}
            }
        \end{kase}
        \item Suppose that $c$ belongs to two triangles in $G[U]$ that are parts of $\D$.
        \begin{kase}
            \item Suppose that $abc, cde \in \D$. Then the subgraph of $I(\D)$ induced by $\sset{abc, cde, \ol{ad}, \ol{be}}$ is isomorphic to the diamond graph, which is the line graph of the paw graph. By \cref{thm:whitney}, the vertices $abc, cde, \ol{ad}, \ol{be}$ of $I(\D)$ represent the edges of a paw graph as follows.
            \figurehc[2]{
                \defpolar{a/0/0,b/0/1,c/150/1,d/210/1}
                \drawedgesbadge{a/b/\ol{be},a/c/cde,a/d/abc,c/d/\ol{ad}}
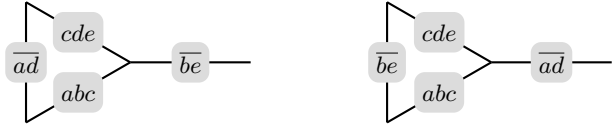
                \begin{scope}[shift={(3,0)}]
                    \defpolar{a/0/0,b/0/1,c/150/1,d/210/1}
                    \drawedgesbadge{a/b/\ol{ad},a/c/cde,a/d/abc,c/d/\ol{be}}
                \end{scope}
            }
            Up to symmetry, assume that $abc, cde, \ol{ad}$ represent the edges of a triangle in $F'$. Since $\ol{ad}$ does not represent a pendant edge of $F'$, it is a triangle part of $\D$. Let $f$ be a vertex of $G$ such that $\ol{ad} = adf$. Since the diamond in $F$ is induced, it follows that $a,c,d,f$ share a common endpoint, and hence $cf \in E(G)$.
            \figurehc[1.4]{
                \defpolar{1/0/0,2/30/1,3/90/1,4/150/1,5/90/2}
                \drawlabeledges{1/2/e,1/3/c,1/4/b,2/3/d,3/4/a,3/5/f}
            }
            The part $\ol{cf}$ then intersects with $abc, cde, adf$, so these four vertices of $I(\D)$ represent the edges of a star of size $4$ in $F'$, contradicting the fact that $abc, cde, adf$ form a triangle in $F'$.
            \item Suppose that $acd, bce \in \D$. Since $c \notin R$, it belongs to $3$ parts of $\D$. In particular, $c$ is adjacent to another vertex $f$ in $G$. One can check that $\ol{ab}, \ol{de}, \ol{cf}$ are pairwise disjoint.
            \figurehc[1.2]{
                \defpolar{0/0/0,1/0/1,2/60/1,3/-60/1}
                \drawedges{0/1,0/2,0/3}
                \nodesvertexbadge{0/acd,1/\ol{de},2/\ol{ab},3/\ol{cf}}
            }
            Then the subgraph of $I(\D)$ induced by $\sset{acd, \ol{ab}, \ol{de}, \ol{cf}}$ is isomorphic to \cref{fig:beineke-a}, which contradicts \cref{thm:beineke}. \qedhere
        \end{kase}
    \end{kase}
\end{proof}

Based on the clique number of $F$, we analyze three cases.

\begin{lemma} \label{lem:clique-4}
    Let $F$ be a connected graph and let $R$ be a nonempty subset of pendant edges of $F$. Suppose that the rooted graph $L(F)_R$ admits a valid decomposition $\D$, and there exists a graph $F'$ such that $I(\D) = L(F')$ and each edge part of $\D$, viewed as a vertex of $I(\D)$, represents a pendant edge of $F'$. If the diamond graph is not an induced subgraph of $F$, and the clique number of $F$ is at least $4$, then the incidence graph of $\D$ is isomorphic to \cref{fig:g6-d}.
\end{lemma}

\begin{proof}
    Set $G = L(F)$. Suppose for a moment that $F$ contains a clique $K_5$ of order $5$. Since $F$ is connected, and $F$ has at least one pendant edge, there exists an edge share a single vertex, say $u$, with $K_5$ in $F$. Let $a$ be an edge of $K_5$ with $u$ as an endpoint. We have $\deg_G(a) \ge 7$. However, $a$ belongs to at most $3$ parts of $\D$.

    Hereafter, we assume that the clique number of $F$ is $4$. Suppose that $F$ contains a clique $K_4$ of order $4$ whose vertices are $u,v,w,x$ and whose edges are $a,b,c,d,e,f$ as shown below.
    \figurehc[1.4]{
        \defpolar{1/0/0,2/90/1,3/210/1,4/-30/1}
        \drawlabeledges{1/2/a,1/3/f,1/4/e,3/4/d,4/2/c,2/3/b}
        \nodeslabelvertex{1/u,2/v,3/w,4/x}
    }
    Set $U = \sset{a,b,c,d,e,f}$. Let $\D'$ be the set of triangles in $G[U]$ that are parts of $\D$. Since the degree of each vertex of $U$ is at least $4$ in $G$, each belongs to a triangle part of $\D'$, which implies that $\D'$ contains at least two triangle parts. Since $G[U]$ contains $12$ edges, we know that $\D'$ contains at most $4$ triangle parts.

    \begin{kase}
        \item Suppose that $\D'$ contains exactly $2$ triangle parts. Since each vertex of $U$ belongs to a triangle part of $\D'$, up to symmetry, we may assume that $\D' = \sset{abc, def}$. Then the subgraph of $I(\D)$ induced by $\sset{abc, \ol{af}, \ol{bd}, \ol{ce}}$ is isomorphic to \cref{fig:beineke-a}, which contradicts \cref{thm:beineke}.
        \figurehc[1.2]{
            \defpolar{0/0/0,1/0/1,2/60/1,3/-60/1}
            \drawedges{0/1,0/2,0/3}
            \nodesvertexbadge{0/abc,1/\ol{bd},2/\ol{af},3/\ol{ce}}
        }
        \item Suppose that $\D'$ contains exactly $3$ triangle parts. Without loss of generality, we may assume that $abc \in \D'$ or $def \in \D'$. Note that we cannot have both $abc \in \D'$ and $def \in \D'$ because, if so, the third triangle part in $\D'$ would have to intersect with $abc$ or $def$ at two vertices.
        \begin{kase}
            \item Suppose that $abc \in \D'$ and $def \notin \D'$. Each of the other two triangle parts in $\D'$ has one vertex from $abc$ and two vertices from $def$. Up to symmetry, we may assume that $D' = \sset{abc, bdf, cde}$, and so $\ol{ae}, \ol{af}, \ol{ef}$ are three distinct parts of $\D$. One can check that the subgraph of $I(\D)$ induced by $\sset{abc, bdf, cde, \ol{ae}, \ol{af}, \ol{ef}}$ is isomorphic to the line graph of the complete graph $K_4$ of order $4$. By \cref{thm:whitney}, the vertices $abc, bdf, cde, \ol{ae}, \ol{af}, \ol{ef}$ of $I(\D)$ represent the edges of a complete graph $K_4$ as follows.
            \figurehc[2]{
                \defpolar{1/0/0,2/90/1,3/210/1,4/-30/1}
                \drawedgesbadge{1/2/abc,1/3/\ol{af},1/4/\ol{ae},3/4/\ol{ef},4/2/cde,2/3/bdf}
            }
            Thus $\ol{ef}$ does not represent a pendant edge in $F'$, and so $\ol{ef}$ is a triangle part of $\D$. Let $g$ be an edge of $F$ such that $\ol{ef} = efg$. Clearly $ag \in E(G)$. Since $a$ already belongs to three parts $abc, \ol{ae}, \ol{af}$ of $\D$, it must be the case that either $\ol{ae} = aeg$ or $\ol{af} = afg$, which contradicts $efg \in \D$.

            \item Suppose that $abc \notin \D'$ and $def \in \D'$. Each of the other two triangle parts in $\D'$ has one vertex from $def$ and two vertices from $abc$. Up to symmetry, we may assume that $\D' = \sset{abf, ace, def}$, and so $\ol{bc}, \ol{bd}, \ol{cd}$ are three distinct parts of $\D$. One can check that the subgraph of $I(\D)$ induced by $\sset{abf, ace, def, \ol{bc}, \ol{cd}, \ol{bd}}$ is isomorphic to the line graph of the complete graph $K_4$ of order $4$. By \cref{thm:whitney}, the vertices $abf, ace, def, \ol{bc}, \ol{cd}, \ol{bd}$ of $I(\D)$ represent the edges of a complete graph $K_4$ as follows.
            \figurehc[2]{
                \defpolar{1/0/0,2/90/1,3/210/1,4/-30/1}
                \drawedgesbadge{1/2/abf,1/3/\ol{bc},1/4/\ol{bd},3/4/\ol{cd},4/2/def,2/3/ace}
            }
            Thus none of $\ol{bc}, \ol{cd}, \ol{bd}$ represents a pendant edge in $F'$, and so they are triangle parts of $\D$. Let $g, h, i$ be edges of $F$ such that $\ol{bc} = bcg, \ol{cd} = cdh, \ol{bd} = bdi$. These edges are incident with $v,x,w$, respectively, and their other endpoints lie outside $\sset{u,v,w,x}$. Write $g = vy$. Suppose for a moment that there exists another edge $j$ incident with $u$ in $F$. Since $a$ already belongs to $abf, ace, \ol{ag}$, we have $\ol{ag} = agj$, and so $j = uy$. Since the diamond graph is not an induced subgraph of $F$, we have $yw, yx \in E(F)$. Then $\sset{u,v,w,x,y}$ is a clique of order $5$ in $F$, which contradicts the assumption that the clique number of $F$ is $4$. Thus $a,e,f$ are the only edges incident with $u$ in $F$. Since the diamond graph is not an induced subgraph of $F$, the edges $g, h, i$ form a matching in $F$ as shown below.
            \figurehc[1.4]{
                \defpolar{1/0/0,2/90/1,3/210/1,4/-30/1,5/90/2,6/210/2,7/-30/2}
                \drawlabeledges{1/2/a,1/3/f,1/4/e,3/4/d,4/2/c,2/3/b,2/5/g,3/6/i,4/7/h}
                \nodeslabelvertex{1/u,2/v,3/w,4/x,5/y}
            }
            Since $b$ already belongs to $abf, bcg, bdi$, we conclude that $a,b,c,g$ are the only edges incident with $v$, and so $\ol{ag} = ag$. Similarly, $eh, fi \in \D$.

            \begin{kase}
                \item Suppose that not all of $g,h,i$ are pendant edges of $F$. Up to symmetry, assume that $g$ is not a pendant edge. Let $j$ be another edge of $F$ incident with $y$. Then $g$ belongs to three parts $ag, bcg, \ol{gj}$ of $\D$. The subgraph of $I(\D)$ induced by $\sset{abf, bcg, bdi, ag, fi, \ol{gj}}$ is isomorphic to \cref{fig:beineke-h}, which contradicts \cref{thm:beineke}.
                \figurehc[1.2]{
                    \defpolar{0/0/0,1/60/1,2/120/1,3/180/1,4/240/1}
                    \coordinate (5) at (210:{sqrt(3)});
                    \drawedges{0/2,1/2,0/3,0/4,4/5,2/3,3/4,0/1,3/5}
                    \nodesvertexbadge{0/bcg,1/\ol{gj},2/ag,3/abf,4/bdi,5/fi}
                }
                \item Suppose that $g, h, i$ are pendant edges of $F$. Since $u,v,w,x$ have no further incident edges and the other endpoints of $g,h,i$ are leaves, connectedness implies that $\sset{a,\dots, i}$ are all the edges of $F$. Hence $\D = \sset{abf, ace, def, bcg, cdh, bdi, ag, eh, fi}$, and its incidence graph is isomorphic to \cref{fig:g6-d}.
                \figurehc[1.2]{
                    \coordinate (A) at (0,0);
                    \coordinate (B) at (0,1);
                    \coordinate (C) at ($(A)+(150:1)$);
                    \coordinate (D) at ($(C)+(90:2)$);
                    \coordinate (E) at ($(B)+(150:2)$);
                    \coordinate (F) at ($(E)+(90:1)$);
                    \coordinate (G) at ($(E)+(-120:.866)$);
                    \coordinate (H) at ($(C)+(180:.866)$);
                    \coordinate (I) at ($(A)+(30:1)$);
                    \coordinate (J) at ($(I)+(90:2)$);
                    \coordinate (K) at ($(B)+(30:2)$);
                    \coordinate (L) at ($(K)+(90:1)$);
                    \coordinate (M) at ($(K)+(-60:.866)$);
                    \coordinate (N) at ($(I)+(0:.866)$);
                    \coordinate (O) at ($(D)+(-30:1)$);
                    \coordinate (R) at ($(F)+(30:2)$);
                    \coordinate (P) at ($(O)!1/3!(R)$);
                    \coordinate (Q) at ($(R)!1/3!(O)$);
                    \draw (E)--(F)--(D)--(C)--(A)--(B)--(E)--(G)--(H)--(C);
                    \draw (A)--(I)--(J)--(L)--(K)--(B);
                    \draw (I)--(N)--(M)--(K);
                    \draw (F)--(R)--(L);
                    \draw (D)--(O)--(J);
                    \draw (O)--(P)--(Q)--(R);
                    \nodesvertexbadge{B/def,C/bdi,F/abf,G/fi,I/cdh,L/ace,M/eh,O/bcg,Q/ag}
                    \nodesvertexbadgehollow{A/d,D/b,E/f,H/i,J/c,K/e,N/h,P/g,R/a}
                }
            \end{kase}
        \end{kase}
        \item Suppose that $\D'$ contains exactly $4$ triangle parts. There are only two different triangle decomposition of $G[U]$, that is $\D' = \sset{abc, aef, bdf, cde}$ or $\D' = \sset{abf, ace, bcd, def}$. In either case, every vertex of $U$ belongs to two triangle parts in $\D'$. Since $a \notin R$, it belongs to $3$ parts of $\D$. Since $a$ only belongs to $2$ parts in $\D'$, there exists another edge, say $g$, of $F$ such that $ag \in E(G)$. Up to symmetry, we may assume that $g$ is incident with $v$. Let $y$ be the vertex of $F$ such that $g = vy$.
        \figurehc[1.2]{
            \defpolar{1/0/0,2/90/1,3/210/1,4/-30/1,5/90/2}
            \drawlabeledges{1/2/a,1/3/f,1/4/e,3/4/d,4/2/c,2/3/b,2/5/g}
            \nodeslabelvertex{1/u,2/v,3/w,4/x, 5/y}
        }
        Since the clique number of $F$ is $4$, the vertex $y$ cannot be adjacent to all vertices in $\sset{u,w,x}$ in $F$. Furthermore, since the diamond graph is not an induced subgraph of $F$, one can deduce that $y$ is not adjacent to any of the vertices in $\sset{u,w,x}$ in $F$. Since $f \notin R$ and $f$ only belongs to $2$ parts in $\D'$, there exists another edge, say $h$, of $F$ such that $fh \in E(G)$. Since $yu, yw \notin E(F)$, we know that $h$ is not incident with $y$ in $F$.
        \begin{kase}
            \item If $h$ is incident with $u$ in $F$, then $a$ belongs to $2$ parts of $\D'$ and $2$ parts $\ol{ag}, \ol{ah}$ of $\D$, which is a contradiction.
            \item If $h$ is incident with $w$ in $F$, then $b$ belongs to $2$ parts of $\D'$ and $2$ parts $\ol{bg}, \ol{bh}$ of $\D$, which is a contradiction. \qedhere
        \end{kase}
    \end{kase}
\end{proof}

\begin{lemma} \label{lem:clique-3}
    Let $F$ be a connected graph and let $R$ be a subset of pendant edges of $F$. Suppose that the rooted graph $L(F)_R$ admits a valid decomposition $\D$, and there exists a graph $F'$ such that $I(\D) = L(F')$ and each edge part of $\D$, viewed as a vertex of $I(\D)$, represents a pendant edge of $F'$. If the diamond graph is not an induced subgraph of either $F$ or $F'$, and the clique number of $F$ is $3$, then the incidence graph of $\D$ is isomorphic to \cref{fig:g6-c}.
\end{lemma}

\begin{proof}
    Set $G = L(F)$. Suppose that $F$ contains a triangle, whose vertices are $u,v,w$, and whose edges are $a,b,c$ as shown below.
    \figurehc{
        \defpolar{2/90/1,3/210/1,4/-30/1}
        \drawlabeledges{3/4/c,4/2/a,2/3/b}
        \nodeslabelvertex{3/u,4/v,2/w}
    }
    Since the clique number of $F$ is $3$, and the diamond graph is not an induced subgraph of $F$, every pair of edges that are incident with $u$, $v$, or $w$ cannot share endpoints outside $\sset{u,v,w}$. Up to symmetry, we may assume that $\deg_F(u) \le \deg_F(v) \le \deg_F(w)$.
    \begin{kase}
        \item Suppose that $\deg_F(u) = \deg_F(v) = 2$. Since $c \notin R$, it belongs to $3$ parts of $\D$. However, $ac$ and $bc$ are the only edges of $G$ that are incident with $c$.

        \item Suppose that $\deg_F(u) = 2$ and $\deg_F(v) = \deg_F(w) = 3$. Let $d$ and $e$ be the other edges of $F$ incident with $v$ and $w$ respectively.
        \figurehc{
            \defpolar{2/90/1,3/210/1,4/-30/1}
            \coordinate (5) at ($(2) + (0:{sqrt(3)})$);
            \coordinate (7) at ($(4) + (0:{sqrt(3)})$);
            \drawlabeledges{3/4/c,4/2/a,2/3/b,2/5/e,4/7/d}
            \nodeslabelvertex{3/u,4/v,2/w}
        }
        Since $b$ and $c$ each belongs to $3$ parts of $\D$, it must be the case that $ab, bc, be, ac, cd \in \D$. Now $a$ belongs to four distinct parts $ab, ac, \ol{ad}, \ol{ae}$ of $\D$.

        \item Suppose that $\deg_F(u) = 2$, $\deg_F(v) \ge 3$, and $\deg_F(w) \ge 4$. Let $d$ be an edge of $F$ incident with $v$, and let $e, f$ be edges of $F$ incident with $w$.
        \figurehc{
            \defpolar{2/90/1,3/210/1,4/-30/1}
            \coordinate (5) at ($(2) + (180:{sqrt(3)})$);
            \coordinate (6) at ($(2) + (0:{sqrt(3)})$);
            \coordinate (7) at ($(4) + (0:{sqrt(3)})$);
            \drawlabeledges{3/4/c,4/2/a,2/3/b,2/5/f,4/7/d,2/6/e}
            \nodeslabelvertex{3/u,4/v,2/w}
        }
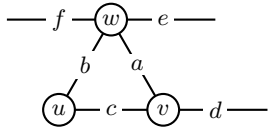
        \begin{kase}
            \item Suppose that $aef \notin \D$. Since $a$ belongs to at most $3$ parts of $\D$, one can deduce that either $abe \in \D$ or $abf \in \D$. Up to symmetry, we may assume that $abe \in \D$. Thus $a$ belongs to $3$ parts $abe, \ol{af}, acd$ of $\D$, and moreover, $\ol{bf} \in \D$. Since $abe, acd \in \D$, and $b, c$ are the only edges incident with $u$ in $F$, we have $bc \in \D$. The subgraph of $I(\D)$ induced by $\sset{abe, \ol{af}, acd, \ol{bf}, bc}$ is isomorphic to the line graph of the diamond graph. Since the diamond graph cannot be an induced subgraph of $F'$, the subgraph of $F'$ induced by the vertices of the diamond graph has to be the complete graph $K_4$. As a consequence, there exists a part $\sigma$ of $\D$ that is adjacent to $\ol{af}, acd, \ol{bf}, bc$ but not $abe$ in $I(\D)$.
            \figurehc[2]{
                \defpolar{0/0/0,1/90/1,2/210/1,3/-30/1}
                \drawedgesbadge{0/1/\sigma,0/2/acd,0/3/bc,1/2/\ol{af},2/3/abe,3/1/\ol{bf}}
            }
            Since $\sigma$ and $abe$ are disjoint but $\sigma$ and $bc$ intersect, we know that $b \notin \sigma$ but $c \in \sigma$. Since $\sigma \neq acd$, there exists $g \in \sigma$ such that $cg \in E(G)$. Since $b$ and $c$ are the only edges of $F$ incident with $u$, we know that $g$ is incident with $v$ in $F$, and so $ag \in E(G)$. Now $a$ belongs to $4$ parts $abe, \ol{af}, acd, \ol{ag}$ of $\D$.
            \item Suppose that $aef \in \D$. Then $\ol{ab}, \ol{be}, \ol{bf}$ are three distinct parts of $\D$. Since $b$ belongs to at most $3$ parts, it must be case that $\ol{ab} = abc$, and so $\ol{cd}$ and $\ol{ad}$ are two distinct parts of $\D$. Since $c \notin R$, it belongs to $3$ parts of $\D$, and so there exists another edge $g$ incident with $v$ in $F$.
            \figurehc{
                \defpolar{2/90/1,3/210/1,4/-30/1}
                \coordinate (5) at ($(2) + (180:{sqrt(3)})$);
                \coordinate (6) at ($(2) + (0:{sqrt(3)})$);
                \coordinate (7) at ($(4) + (-30:{sqrt(3)})$);
                \coordinate (8) at ($(4) + (30:{sqrt(3)})$);
                \drawlabeledges{3/4/c,4/2/a,2/3/b,2/5/f,4/7/d,2/6/e,4/8/g}
                \nodeslabelvertex{3/u,4/v,2/w}
            }
            Since $a$ already belongs to two parts $aef, abc$ of $\D$, the third part $\ol{ad} = adg$. Therefore $\deg_F(v) = \deg_F(w) = 4$, and so $\ol{be} = be, \ol{bf} = bf, \ol{cd} = cd, \ol{cg} = cg$.
            \begin{kase}
                \item Suppose that $F$ consists of only $a,b,c,d,e,f,g$ as its edges. We conclude that $\D = \sset{abc, adg, aef, be, bf, cd, cg}$. Then the incidence graph of $\D$ is isomorphic to \cref{fig:g6-c}.
                \figurehc{
                    \coordinate (A) at (0,0);
                    \coordinate (B) at (0,1);
                    \coordinate (C) at ($(A)+(150:1)$);
                    \coordinate (D) at ($(C)+(90:2)$);
                    \coordinate (E) at ($(B)+(150:2)$);
                    \coordinate (F) at ($(E)+(90:1)$);
                    \coordinate (G) at ($(E)+(-120:.866)$);
                    \coordinate (H) at ($(C)+(180:.866)$);
                    \coordinate (I) at ($(A)+(30:1)$);
                    \coordinate (J) at ($(I)+(90:2)$);
                    \coordinate (K) at ($(B)+(30:2)$);
                    \coordinate (L) at ($(K)+(90:1)$);
                    \coordinate (M) at ($(K)+(-60:.866)$);
                    \coordinate (N) at ($(I)+(0:.866)$);
                    \draw (E)--(F)--(D)--(C)--(A)--(B)--(E)--(G)--(H)--(C);
                    \draw (A)--(I)--(J)--(L)--(K)--(B);
                    \draw (I)--(N)--(M)--(K);
                    \nodesvertexbadge{B/abc,C/adg,F/cd,G/cg,I/aef,L/be,M/bf}
                    \nodesvertexbadgehollow{A/a,D/d,E/c,H/g,J/e,K/b,N/f}
                }
                \item Suppose there exists another edge, say $h$, of $F$. Since $\deg_F(u) = 2$, $\deg_F(v) = \deg_F(w) = 4$, and $F$ is connected, we may assume that $h$ is adjacent to $d$, $e$, $f$, or $g$ in $G$. Up to symmetry, we may assume that $dh \in E(G)$. Let $x$ be the vertex of $F$ such that $d = vx$. Then $h$ is incident to $x$ in $F$. Since $adg \in \D$, we have $\ol{dh} \neq dgh$. Then the subgraph of $I(\D)$ induced by $\sset{\ol{dh}, cd, adg, abc, aef, be}$ is isomorphic to \cref{fig:beineke-h}, which contradicts \cref{thm:beineke}.
                \figurehc[1.2]{
                    \defpolar{0/0/0,1/60/1,2/120/1,3/180/1,4/240/1}
                    \coordinate (5) at (210:{sqrt(3)});
                    \drawedges{0/2,1/2,0/3,0/4,4/5,2/3,3/4,0/1,3/5}
                    \nodesvertexbadge{0/adg,1/\ol{dh},2/cd,3/abc,4/aef,5/be}
                }
            \end{kase}
        \end{kase}

        \item Suppose that $\deg_F(u) = \deg_F(v) = \deg_F(w) = 3$. Let $d, e, f$ be the other edges incident with $u,v,w$ respectively.
        \figurehc{
            \defpolar{2/90/1,3/210/1,4/-30/1}
            \coordinate (5) at ($(3) + (180:{sqrt(3)})$);
            \coordinate (6) at ($(2) + (180:{sqrt(3)})$);
            \coordinate (7) at ($(4) + (0:{sqrt(3)})$);
            \drawlabeledges{3/4/c,4/2/a,2/3/b,3/5/d,4/7/e,2/6/f}
            \nodeslabelvertex{3/u,4/v,2/w}
        }
        \begin{kase}
            \item If $abc \in \D$, then $\ol{ae}, \ol{bf}, \ol{cd}$ are three distinct parts of $\D$, and so the subgraph of $I(\D)$ induced by $\sset{abc, \ol{ae}, \ol{bf}, \ol{cd}}$ is isomorphic to \cref{fig:beineke-a}, which contradicts \cref{thm:beineke}.
            \figurehc[1.2]{
                \defpolar{0/0/0,1/0/1,2/60/1,3/-60/1}
                \drawedges{0/1,0/2,0/3}
                \nodesvertexbadge{0/abc,1/\ol{bf},2/\ol{ae},3/\ol{cd}}
            }
            \item Suppose that $abc \notin \D$. Then $\ol{ab}, \ol{ac}, \ol{bc}$ are three distinct parts of $\D$. Since $a \notin R$, it belongs to $3$ parts of $\D$, and so either $abf \in \D$ or $ace \in \D$ exclusively. Similarly, either $abf \in \D$ or $bcd \in \D$ exclusively, and either  $ace \in \D$ or $bcd \in \D$ exclusively. However, it is impossible to select $abf, ace, bcd$ to $\D$ to satisfy these constraints.
        \end{kase}

        \item Suppose that the $\deg_F(u), \deg_F(v) \ge 3$ and $\deg_F(w) \ge 4$. Let $d$ and $e$ be edges of $F$ incident with $u$ and $v$ respectively, and let $f$ and $g$ incident with $w$ in $F$.
        \figurehc{
            \defpolar{2/90/1,3/210/1,4/-30/1}
            \coordinate (5) at ($(3) + (180:{sqrt(3)})$);
            \coordinate (6) at ($(2) + (180:{sqrt(3)})$);
            \coordinate (7) at ($(4) + (0:{sqrt(3)})$);
            \coordinate (8) at ($(2) + (0:{sqrt(3)})$);
            \drawlabeledges{3/4/c,4/2/a,2/3/b,3/5/d,2/6/f,2/8/g,4/7/e}
            \nodeslabelvertex{3/u,4/v,2/w}
        }
        Since $a$ belongs to at most $3$ parts of $\D$, it belongs to a triangle part of $\D$ that is contained in $\sset{a,b,f,g}$. Similarly, so does $b$. It must be the case that $abf \in \D$ or $abg \in \D$. Up to symmetry, we may assume that $abf \in \D$. It follows that $\ol{ag}, \ol{bg}, \ol{fg}$ are three distinct parts of $\D$. Since $a$ already belongs to $abf, \ol{ag} \in \D$, and $ac, ae \in E(G)$, we must have $ace \in \D$. Similarly, $bcd \in \D$. The subgraph of $I(\D)$ induced by $\sset{abf, \ol{ag}, \ol{bg}, ace, bcd}$ is isomorphic to the line graph of the diamond graph.
        \cref{thm:whitney} implies that the diamond graph is a subgraph of $F'$. Since the diamond graph cannot be an induced subgraph of $F'$, the subgraph of $F'$ induced by the vertices of the diamond graph has to be the complete graph $K_4$. As a consequence, there exists a part $\sigma$ of $\D$ that is adjacent to $\ol{ag}, \ol{bg}, ace, bcd$ but not $abf$ in $I(\D)$.
        \figurehc[2]{
            \defpolar{0/0/0,1/90/1,2/210/1,3/-30/1}
            \drawlabeledges{0/1/\sigma,0/2/ace,0/3/bcd,1/2/\ol{ag},2/3/abf,3/1/\ol{bg}}
        }
        Since $\sigma$ and $abf$ are disjoint, we have $a \notin \sigma$. Since $\sigma$ and $ace$ intersect, we have $c \in \sigma$ or $e \in \sigma$. Since $cg \notin E(G)$ and $eg \notin E(G)$, we have $g \notin \sigma$. Since $\sigma$ and $\ol{ag}$ intersect, there exists $h \in \sigma$ such that $\ol{ag} = agh$. Recall that every pair of edges that are incident with $u$, $v$ or $w$ cannot share endpoints outside $\sset{u,v,w}$. We know that $h$ is incident with $w$, and $ch, eh \notin E(G)$, which contradicts $c \in \sigma$ or $e \in \sigma$. \qedhere
    \end{kase}
\end{proof}

\begin{lemma} \label{lem:clique-2}
    Let $F$ and $F'$ be trees, and let $R$ be a subset of pendant edges of $F$. If the rooted graph $L(F)_R$ admits a valid decomposition $\D$ with at least one edge part such that the intersection graph $I(\D)$ is the line graph of $F'$, then the incidence graph of $\D$ is isomorphic to \cref{fig:g6-a} or \cref{fig:g6-b}.
\end{lemma}

\begin{proof}
    Set $G = L(F)$.
    \begin{kase}
        \item Suppose that the maximum degree of $F$ is at least $5$. Let $v$ be a vertex of $F$ of degree at least $5$, and let $a,b,c,d,e$ be five edges incident with $v$ in $F$. Set $U = \sset{a,b,c,d,e}$. Then the subgraph of $G$ induced by $U$ is isomorphic to the complete graph $K_5$. Clearly, each vertex in $U$ belongs to a triangle part of $\D$ that is a subset of $U$. Therefore, there are at least two triangle parts of $\D$ that are subsets of $U$. Without loss of generality, we assume that $abc, ade \in \D$, and so $\ol{bd}, \ol{be}, \ol{cd}, \ol{ce}$ are four distinct parts of $\D$.
        \begin{kase}
            \item If $\ol{bd}$ and $\ol{ce}$ are disjoint or $\ol{be}$ and $\ol{cd}$ are disjoint, then the subgraph of $I(\D)$ induced by $\sset{abc, ade, \ol{bd}, \ol{cd}, \ol{ce}}$ or $\sset{abc, ade, \ol{bd}, \ol{be}, \ol{cd}}$ is isomorphic to \cref{fig:beineke-c}, which contradicts \cref{thm:beineke}.
            \figurehc[1.2]{
                \defpolar{0/0/0,1/0/1,2/120/1,3/240/1,4/180/1}
                \drawedges{0/1,1/2,0/2,0/3,0/4,1/3,2/3,2/4,3/4}
                \nodesvertexbadge{0/ade,1/\ol{ce},2/abc,3/\ol{cd},4/\ol{bd}}
                \begin{scope}[shift={(4,0)}]
                    \defpolar{0/0/0,1/0/1,2/120/1,3/240/1,4/180/1}
                    \drawedges{0/1,1/2,0/2,0/3,0/4,1/3,2/3,2/4,3/4}
                    \nodesvertexbadge{0/abc,1/\ol{cd},2/ade,3/\ol{bd},4/\ol{be}}
                \end{scope}
            }
            \item Suppose that $\ol{bd}$ and $\ol{ce}$ intersect and $\ol{be}$ and $\ol{cd}$ intersect. There exist edges $f$ and $g$ of $F$ such that $\ol{bd} = bdf, \ol{ce} = cef, \ol{be} = beg, \ol{cd} = cdg$, and so both $f$ and $g$ are incident to $v$ in $F$. Since $f$ already belongs to two parts $bdf$ and $cef$ of $\D$, we have $afg \in \D$. Since each of $a,b,c,d,e,f,g$ already belongs to three triangle parts of $\D$, no other vertex of $G$ is adjacent to any of them. Since $F$ is connected, it consists of only these seven edges, and all parts of $\D$ are triangle parts, which contradicts the condition that $\D$ has at least one edge part.
        \end{kase}
        \item Suppose that the maximum degree of $F$ is $4$. Let $v$ be a vertex of $F$ of degree $4$, and let $a,b,c,d$ be the four edges incident with $v$ in $F$.
        \begin{kase}
            \item Suppose that the graph $F$ consists of only $a,b,c,d$ as its edges.
            \figurehc[1.4]{
                \defpolar{0/0/0,1/0/1,2/60/1,3/120/1,4/0/-1}
                \drawlabeledges{0/1/a,0/2/b,0/3/c,0/4/d}
                \nodeslabelvertex{0/v}
            }
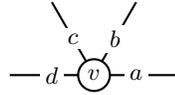
            \begin{kase}
                \item Suppose that $\D$ consists of edge parts only, that is, $\D = \sset{ab,ac,ad,bc,bd,cd}$. Thus $I(\D)$ is the line graph of the complete graph $K_4$. By \cref{thm:whitney}, the graph $F'$ is isomorphic to $K_4$, which contradicts the assumption that $F'$ is a tree.
                \item If $\D$ consists of at least one triangle part, then, up to symmetry, we may assume that $abc \in \D$, and so $\D = \sset{abc, ad, bd, cd}$. Then the incidence graph of $\D$ is isomorphic to \cref{fig:g6-b}.
                \figurehc{
                    \coordinate (A) at (0,0);
                    \coordinate (B) at ($(A)+(30:1)$);
                    \coordinate (C) at ($(A)+(-30:1)$);
                    \coordinate (D) at ($(C)+(30:2)$);
                    \coordinate (E) at ($(B)+(-30:2)$);
                    \coordinate (F) at ($(E)+(30:1)$);
                    \coordinate (B1) at ($(B)+(60:.866)$);
                    \coordinate (D1) at ($(D)+(120:.866)$);
                    \draw (B)--(E)--(F)--(D)--(C)--(A)--(B)--(B1)--(D1)--(D);
                    \nodesvertexbadge{B/abc,C/ad,F/bd,D1/cd}
                    \nodesvertexbadgehollow{A/a,D/d,E/b,B1/c}
                }
            \end{kase}
            \item Suppose that there exists another edge, say $e$, of $F$ adjacent to one of $a,b,c,d$ in $G$. Up to symmetry, we may assume that $de \in E(G)$. Let $u$ be a vertex of $F$ such that $d = uv$. Thus $e$ is incident with $u$ in $F$.
            \figurehc[1.4]{
                \defpolar{0/0/0,1/0/1,2/60/1,3/120/1,4/0/-1}
                \coordinate (5) at ($(3)+(-1,0)$);
                \drawlabeledges{0/1/a,0/2/b,0/3/c,0/4/d,4/5/e}
                \nodeslabelvertex{0/v,4/u}
            }
            Since $d$ belongs to at most $3$ parts of $\D$, and it already belongs to $\ol{de}$, it must belong to a triangle part of $\D$ that is a subset of $\sset{a,b,c,d}$. Up to symmetry, we may assume that $bcd \in \D$, and so $\ol{ab}, \ol{ac}, \ol{ad}$ are $3$ distinct parts of $\D$. Since $F$ is a tree, and its maximum degree is $4$, we know that $\ol{ab} = ab, \ol{ac} = ac, \ol{ad} = ad$, and they are disjoint from $\ol{de}$. Therefore, the subgraph of $I(\D)$ induced by $\sset{\ol{de}, bcd, ab, ad}$ is isomorphic to the diamond graph, which is impossible since two vertices in the line graph of the tree $F'$ that are not adjacent have at most one common neighbor.
        \end{kase}
        \item Suppose that the maximum degree of $F$ is at most $3$. Let $a$ be a pendant edge of $F$. Since $a$ belongs to at least $2$ parts of $\D$, one endpoint of $a$, say $v$, is of degree exactly $3$. Let $b$ and $c$ be the other two edges incident with $v$ in $F$.
        \begin{kase}
            \item Suppose that the graph $F$ consists of only $a,b,c$ as its edges. Since $a$ belongs to at least $2$ parts, we have $\D = \sset{ab,ac,bc}$. Then the incidence graph of $\D$ is isomorphic to \cref{fig:g6-a}.
            \figurehc{
                \coordinate (A) at (0,0);
                \coordinate (B) at ($(A)+(30:1)$);
                \coordinate (C) at ($(A)+(-30:1)$);
                \coordinate (D) at ($(C)+(30:2)$);
                \coordinate (E) at ($(B)+(-30:2)$);
                \coordinate (F) at ($(E)+(30:1)$);
                \draw (A)--(B)--(E)--(F)--(D)--(C)--cycle;
                \nodesvertexbadge{B/ac,C/ab,F/bc}
                \nodesvertexbadgehollow{A/a,D/b,E/c}
            }
            \item Suppose that there exists another edge, say $d$, of $F$, satisfying $bd \in E(G)$ or $cd \in E(G)$. Up to symmetry, we may assume that $cd \in E(G)$. Let $u$ be the vertex of $F$ such that $c = uv$. Thus $d$ is incident with $u$ in $F$.
            \figurehc[1.4]{
                \defxy{0/0/0,1/1/0,2/0/1,3/-1/0,4/-1/1}
                \drawlabeledges{0/1/a,0/2/b,0/3/c,3/4/d}
                \nodeslabelvertex{0/v,3/u}
            }
            Since $a$ is a pendant edge, $\deg_F(v) = 3$, and $a$ belongs to at least $2$ parts of $\D$, we know that $ab, ac, bc \in \D$, and $ab$ and $\ol{cd}$ are disjoint. The subgraph of $I(\D)$ induced by $\sset{ab, ac, bc, \ol{cd}}$ is isomorphic to the diamond graph, which is impossible since two vertices in the line graph of the tree $F'$ that are not adjacent have at most one common neighbor. \qedhere
        \end{kase}
    \end{kase}
\end{proof}

Finally, we assemble all the previous lemmas in this section to classify for the bipartite graphs of girth at least $6$.

\begin{proof}[Proof of \cref{thm:bipartite-girth-6}]
    Suppose that $H$ is a connected subcubic bipartite graph of minimum degree $2$ and girth at least $6$ such that neither \cref{fig:g4-f} nor \cref{fig:g4-g} is a rooted distance-two induced subgraph of $H$, and every rooted distance-two component $G_R$ with $R \neq \varnothing$ is $L(F)$ for some connected graph with petals $F$, with each root in $R$ representing a pendant edge of $F$.

    Since the minimum degree of $H$ is $2$, one of its two rooted distance-two components, say $G_R$, has $R \neq \varnothing$; denote the other by $G_S'$. Write $G = L(F)$ for a connected graph with petals $F$ such that each vertex in $R$ represents a pendant edge of $F$. According to \cref{lem:valid-decomposition}, there exists a valid decomposition $\D$ of $G_R$ such that $G_S'$ and $\iddtwo$ are isomorphic, and $H$ and the incidence graph of $\D$ are isomorphic.

    Since neither \cref{fig:g4-f} nor \cref{fig:g4-g} is a rooted distance-two induced subgraph of $H$, neither is an induced rooted subgraph of $G'_S$. Since $G'_S$ and $\iddtwo$ are isomorphic, \cref{lem:no-petal} implies that $F$ has no petal, and hence $F$ is a simple graph. By \cref{lem:has-edge-parts}, we know that $\D_{\sharp 2} \neq \varnothing$, and so $S \neq \varnothing$.

    Since $S \neq \varnothing$, there exists a connected graph with petals $F'$ such that $G' = L(F')$ and each vertex in $S$ represents a pendant edge of $F'$. A symmetric application of \cref{lem:no-petal} shows that $F'$ is simple. By \cref{lem:chordal,lem:diamond-free}, applied symmetrically, both $F$ and $F'$ are chordal and neither contains the diamond graph as an induced subgraph.

    \begin{kase}
        \item If the clique number of $F$ or $F'$ is at least $4$, then $H$ is isomorphic to \cref{fig:g6-d} by \cref{lem:clique-4}.
        \item If the clique number of $F$ or $F'$ is $3$, then $H$ is isomorphic to \cref{fig:g6-c} by \cref{lem:clique-3}.
        \item Suppose that the clique number of both $F$ and $F'$ is $2$. Since $F$ and $F'$ are connected and chordal, we conclude that both are trees. Therefore $H$ is isomorphic to \cref{fig:g6-a} or \cref{fig:g6-b} by \cref{lem:clique-2}. \qedhere
    \end{kase}
\end{proof}

\section{Non-bipartite graphs} \label{sec:non-bipartite}

\begin{lemma}[Imrich and Pisanski \cite{IP08}] \label{lem:involution}
    Let $H$ be a connected non-bipartite graph. If $H'$ is the bipartite double of $H$, then $H'$ admits an involution $\sigma$ exchanging the two parts of $H'$ such that $v\sigma(v)$ is not an edge of $H'$ for every vertex $v$ of $H'$, and $H'/\sigma$ is isomorphic to $H$. \qed
\end{lemma}

\begin{proof}[Proof of \cref{thm:non-bipartite}]
    Let $H'$ be the bipartite double of $H$. Since $H'$ is connected, $H$ is connected and non-bipartite. According to \cref{lem:involution}, there exists an involution $\sigma$ of $H'$ exchanging the two parts of $H'$ such that $v \sigma(v) \notin E(H')$ for every $v \in V(H')$, and $H'/\sigma$ is isomorphic to $H$.

    \begin{kase}
        \item Suppose that $H'$ is isomorphic to $\hjn$ for some $n \ge 1$. We label the vertices as follows.
        \figurehc[1.2]{
            \defxy{a1/0/0, b1/0/1, c1/1/0, d1/1/1, a2/2/0, b2/2/1, c2/3/0, d2/3/1, a3/4/0, b3/4/1, c3/6/0, d3/6/1, a4/7/0, b4/7/1, c4/8/0, d4/8/1, a5/9/0, b5/9/1, c5/10/0, d5/10/1}
            \defxy{capL0/-1/0, capL1/-1/1, capLc/-2/0, capR0/11/0, capR1/11/1, capRc/12/0}
            \drawedges{a1/c1, c1/b1, b1/d1, d1/a1, a2/c2, c2/b2, b2/d2, d2/a2, a4/c4, c4/b4, b4/d4, d4/a4, a5/c5, c5/b5, b5/d5, d5/a5, c1/a2, d1/b2, c2/a3, d2/b3, c3/a4, d3/b4, c4/a5, d4/b5}
            \drawedges{a1/capL0, capL0/capLc, capLc/capL1, capL1/b1, c5/capR0, capR0/capRc, capRc/capR1, capR1/d5}
            \nodesvertex{capL0,capL1,capLc,capR0,capR1,capRc}
            \draw[dashed] (a3)--(c3);
            \draw[dashed] (b3)--(d3);
            \node[vertex, label=below:$a_1$] at (a1) {};
            \node[vertex, label=above:$a_1'$] at (b1) {};
            \node[vertex, label=below:$b_1$] at (c1) {};
            \node[vertex, label=above:$b_1'$] at (d1) {};
            \node[vertex, label=below:$a_2$] at (a2) {};
            \node[vertex, label=above:$a_2'$] at (b2) {};
            \node[vertex, label=below:$b_2$] at (c2) {};
            \node[vertex, label=above:$b_2'$] at (d2) {};
            \node[vertex, label=below:$a_3$] at (a3) {};
            \node[vertex, label=above:$a_3'$] at (b3) {};
            \nodesvertex{c1,d1,a2,b2,c2,d2,a3,b3,c3,d3,a4,b4,c4,d4,a5,b5}
            \node[vertex, label=below:$b_{n-2}$] at (c3) {};
            \node[vertex, label=above:$b_{n-2}'$] at (d3) {};
            \node[vertex, label=below:$a_{n-1}$] at (a4) {};
            \node[vertex, label=above:$a_{n-1}'$] at (b4) {};
            \node[vertex, label=below:$b_{n-1}$] at (c4) {};
            \node[vertex, label=above:$b_{n-1}'$] at (d4) {};
            \node[vertex, label=below:$a_n$] at (a5) {};
            \node[vertex, label=above:$a_n'$] at (b5) {};
            \node[vertex, label=below:$b_n$] at (c5) {};
            \node[vertex, label=above:$b_n'$] at (d5) {};
        }
        Let $A,B$ be the bipartition classes of $H'$. Their vertices of degree three are respectively $\sset{a_1, a_1', a_2, a_2', \dots, a_n, a_n'}$ and $\sset{b_1, b_1', \dots, b_n, b_n'}$. The only vertices of $H'$ that are of degree three and adjacent to a vertex of degree two are $a_1,a_1',b_n,b_n'$. Since $\sigma$ preserves degrees and exchanges $A,B$, it exchanges $\sset{a_1,a_1'}$ and $\sset{b_n,b_n'}$. Considering the distances of the vertices of degree three to these two sets in $H'$, we then deduce that $\sigma$ exchanges $\sset{b_1,b_1'}$ and $\sset{a_n,a_n'}$, it exchanges $\sset{a_2,a_2'}$ and $\sset{b_{n-1},b_{n-1}'}$, and so on.
        \begin{kase}
            \item If $n = 2m-1$ for some $m \in \N^+$, then $\sigma$ exchanges $\sset{a_m,a_m'}$ and $\sset{b_m,b_m'}$, which contradicts $a_m\sigma(a_m) \notin H'$.
            \item If $n = 2m$ for some $m \in \N^+$, then $\sigma(b_m) = a_{m+1}'$, $\sigma(b_m') = a_{m+1}$, and $H'/\sigma$ is isomorphic to $\mathsf{HJ}'_m$ in \cref{fig:non-bipartite}.
        \end{kase}

        \item Suppose that $H'$ is isomorphic to $\hjone$. We label some of the vertices as follows.
        \figurehc[1.2]{
            \defxy{b0/-2/0,c0/-1/0,d0/-1/1,a1/0/0,b1/0/1,c1/1/0,d1/1/1,a2/2/0,b2/2/1,c2/3/0}
            \drawedges{b1/d0,d0/b0,b0/c0,c0/a1,a1/c1,c1/b1,b1/d1,d1/a1,c1/a2,a2/c2,c2/b2,b2/d1}
            \draw (b0) to[bend right=28.2] (c2);
            \node[vertex, label=above:$a_1$] at (b0) {};
            \node[vertex, label=above:$b_1$] at (c0) {};
            \node[vertex, label=above:$b_1'$] at (d0) {};
            \node[vertex, label=above:$a_2$] at (a2) {};
            \node[vertex, label=above:$a_2'$] at (b2) {};
            \node[vertex, label=above:$b_2$] at (c2) {};
            \nodesvertex{a1,b1,c1,d1}
        }
        Since $b_1, b_1', a_2, a_2'$ are the only vertices of $H'$ that are of degree $2$, and $\sigma$ exchanges the two parts of $H'$, the involution maps $\sset{b_1, b_1'}$ to $\sset{a_2, a_2'}$. Since $a_1$ is the unique common neighbor of $b_1$ and $b_1'$, and $b_2$ is unique common neighbor of $a_2$ and $a_2'$, we conclude that $\sigma(a_1) = b_2$, which contradicts $a_1\sigma(a_1) \notin E(H')$.

        \item Suppose that $H'$ is isomorphic to one of the four graphs in \cref{fig:bipartite-girth-6}. Since the bipartite double of any graph always has even number of edges, $H'$ cannot be \cref{fig:g6-b,fig:g6-c}.
        \begin{kase}
            \item If $H'$ is isomorphic to \cref{fig:g6-a}, then $H$ is isomorphic to $K_3$.
            \item Suppose that $H'$ is isomorphic to \cref{fig:g6-d}. We label its vertices of degree $2$ as follows.
            \figurehc[1.2]{
                \coordinate (A) at (0,0);
                \coordinate (B) at (0,1);
                \coordinate (C) at ($(A)+(150:1)$);
                \coordinate (D) at ($(C)+(90:2)$);
                \coordinate (E) at ($(B)+(150:2)$);
                \coordinate (F) at ($(E)+(90:1)$);
                \coordinate (G) at ($(E)+(-120:.866)$);
                \coordinate (H) at ($(C)+(180:.866)$);
                \coordinate (I) at ($(A)+(30:1)$);
                \coordinate (J) at ($(I)+(90:2)$);
                \coordinate (K) at ($(B)+(30:2)$);
                \coordinate (L) at ($(K)+(90:1)$);
                \coordinate (M) at ($(K)+(-60:.866)$);
                \coordinate (N) at ($(I)+(0:.866)$);
                \coordinate (O) at ($(D)+(-30:1)$);
                \coordinate (R) at ($(F)+(30:2)$);
                \coordinate (P) at ($(O)!1/3!(R)$);
                \coordinate (Q) at ($(R)!1/3!(O)$);
                \draw (E)--(F)--(D)--(C)--(A)--(B)--(E)--(G)--(H)--(C);
                \draw (A)--(I)--(J)--(L)--(K)--(B);
                \draw (I)--(N)--(M)--(K);
                \draw (F)--(R)--(L);
                \draw (D)--(O)--(J);
                \draw (O)--(P)--(Q)--(R);
                \node[vertex, label=left:$a_1$] at (Q) {};
                \node[vertex, label=left:$b_1$] at (P) {};
                \node[vertex, label=left:$a_2$] at (G) {};
                \node[vertex, label=left:$b_2$] at (H) {};
                \node[vertex, label=right:$a_3$] at (M) {};
                \node[vertex, label=right:$b_3$] at (N) {};
                \nodesvertex{A,B,C,D,E,F,I,J,K,L,M,N,O,R}
            }
            Since $a_1\sigma(a_1) \notin E(H')$, we know that $\sigma(a_1) \neq b_1$. Since $\sigma$ is an automorphism of $H'$, $\sigma(a_1b_1) = a_2b_2$ or $\sigma(a_1b_1) = a_3b_3$. Up to symmetry, we may assume that $\sigma(a_1b_1) = a_2b_2$. Since $\sigma$ is an involution, $\sigma(a_2b_2) = a_1b_1$, and so $\sigma(a_3b_3) = a_3b_3$, which contradicts $a_3\sigma(a_3) \notin E(H')$. \qedhere
        \end{kase}
    \end{kase}
\end{proof}

\section{Further remarks} \label{sec:remark}

\subsection{Spectral gap intervals}

A closely related notion is defined as follows. Let $\mathcal{C}$ be a class of graphs, and let $a$ and $b$ denote the infimum and supremum, respectively, of all eigenvalues of graphs in $\mathcal{C}$. An interval $I \subseteq [a,b]$ that is a spectral gap set for $\mathcal{C}$ is called a \emph{spectral gap interval} for $\mathcal{C}$. A spectral gap interval $I$ is said to be \emph{maximal} if it is not properly contained in any other spectral gap interval for $\mathcal{C}$.

For the class $\mathcal{C}$ of connected cubic graphs, many maximal spectral gap intervals fail to be maximal when considered as spectral gap sets for $\mathcal{C}$. For example, Koll\'{a}r, Fitzpatrick, Sarnak, and Houck recently observed that $[-3,-2)$ is a maximal spectral gap interval. However, the set $[-3,-2) \cup (0,1)$ is also a spectral gap set for $\mathcal{C}$. Indeed, it follows from \cite[Section~2.1]{KS21} that the line graph of a subdivision of any cubic graph is itself a cubic graph without eigenvalues in $[-3,-2) \cup (0,1)$.

A similar phenomenon occurs for another well-known maximal spectral gap interval. The celebrated Alon--Boppana bound \cite{N91}, together with the existence of cubic Ramanujan graphs \cite{C92}, implies that $(2\sqrt{2}, 3)$ is a maximal spectral gap interval. Nevertheless, the existence of cubic Ramanujan graphs also shows that $[-3, -2\sqrt2) \cup (2\sqrt{2}, 3)$ is a spectral gap set for $\mathcal{C}$.

To the best of our knowledge, the interval $(-1,1)$ is currently the only interval that is known to be a spectral gap set for $\mathcal{C}$.

\begin{problem}
    Besides $(-1,1)$, find another interval that is maximal as a spectral gap set for the class of connected cubic graphs.
\end{problem}

\subsection{\texorpdfstring{Spectral gap set $(-2,0)$}{Spectral gap set (-2,0)}}

Very recently, Guo and Royle classified in \cite{GR27} all connected cubic graphs without eigenvalues in the open interval $(-2,0)$. It is therefore natural to ask for an analogous classification in the subcubic case.

\begin{problem}
    Classify all connected subcubic graphs without eigenvalues in $(-2,0)$.
\end{problem}

Earlier, Koll\'ar and Sarnak showed in \cite{KS21} that $(-2,0)$ is a spectral gap set for cubic graphs. In contrast to the interval $(-1,1)$, which is maximal under inclusion as a spectral gap set, the interval $(-2,0)$ is not maximal (see \cite[Figure~10]{KS21} and \cite[Section~3]{GR27}).

\subsection{Median eigenvalues}

We end by relating the present work to the line of research on median eigenvalues. Let $G$ be a graph on $n$ vertices with eigenvalues $\lambda_1 \ge \lambda_2 \ge \cdots \ge \lambda_n$. Define the \emph{median eigenvalues} of $G$ to be $\lambda_h$ and $\lambda_\ell$, where $h = \lfloor (n+1)/2 \rfloor$ and $\ell = \lceil (n+1)/2 \rceil$.

Mohar showed in \cite{M16} that, with the exception of the Heawood graph, every connected subcubic bipartite graph has median eigenvalues in the closed interval $[-1,1]$. Since the spectrum of a bipartite graph is symmetric about $0$, a bipartite graph has no median eigenvalues in $(-1,1)$ if and only if it has no eigenvalues in $(-1,1)$. Consequently, combining the bipartite graphs appearing in the classification of Guo and Royle in \cite{GR26} with our classification in \cref{thm:main} yields a complete description of all connected subcubic bipartite graphs whose median eigenvalues avoid $(-1,1)$.

More recently, Acharya, Jeter, and Jiang extended Mohar's result in \cite{AJJ26} to all connected subcubic graphs, showing that, with the exception of the Heawood graph, every connected subcubic graph has median eigenvalues in $[-1,1]$.

Finally, observe that any connected subcubic graph of even order whose median eigenvalues are $-1$ and $+1$ has no eigenvalues in $(-1,1)$. Among the connected cubic graphs without eigenvalues in $(-1,1)$ classified in \cite{GR26}, those with median eigenvalues $-1$ and $+1$ are exactly the Guo--Mohar graphs $\gmn$ for $n \ge 2$, the graph obtained from the complete bipartite graph $K_{3,3}$ by truncating two independent vertices, the graph obtained from the Petersen graph by truncating one vertex, and the eight sporadic graphs that are bipartite other than the Heawood graph. Among the connected subcubic graphs that are not cubic and have no eigenvalues in $(-1,1)$, those with median eigenvalues $-1$ and $+1$ are exactly the bipartite graphs $\hjn$ for $n \ge 1$ in \cref{fig:bipartite-girth-4}, and the six sporadic graphs other than the triangle graph $K_3$.

\section*{Acknowledgements} The second author is grateful to Boris Bukh for providing access to OpenAI's model used in the Lean formalization.

\bibliographystyle{plain}
\bibliography{avoid}

\end{document}